\documentclass[12pt]{report}

\scrollmode
\usepackage{graphicx}
\usepackage{amsfonts}
\usepackage{amsmath}
\usepackage{amssymb}
\usepackage{epsfig}
\usepackage{graphicx}
\usepackage{color}
\newcommand{\thesistitle}{Convergence of Coalescing Nonsimple Random Walks to the Brownian Web}
\newcommand{\myname}{Rongfeng Sun}
\newcommand{\advisor}{Charles M. Newman}
\newcommand{\thesismonth}{January}
\newcommand{\thesisyear}{2005}

\newtheorem{teo}{Theorem}[section]

\newtheorem{prop}[teo]{Proposition}
\newtheorem{lem}{Lemma}[section]
\newtheorem{rmk}[teo]{Remark}
\newtheorem{cor}{Corollary}[section]

\newcommand{\beq}{\begin{equation}}
\newcommand{\eeq}{\end{equation}}
\newcommand{\beqn}{\begin{eqnarray}}
\newcommand{\beqnn}{\begin{eqnarray*}}
\newcommand{\eeqn}{\end{eqnarray}}
\newcommand{\eeqnn}{\end{eqnarray*}}
\newcommand{\bprop}{\begin{prop}}
\newcommand{\eprop}{\end{prop}}
\newcommand{\bteo}{\begin{teo}}
\newcommand{\eteo}{\end{teo}}
\newcommand{\bcor}{\begin{cor}}
\newcommand{\ecor}{\end{cor}}
\newcommand{\brm}{\begin{rmk}}
\newcommand{\erm}{\end{rmk}}
\newcommand{\blem}{\begin{lem}}
\newcommand{\elem}{\end{lem}}

\newcommand{\fr}{\frac}
\renewcommand{\r}{{\mathbb R}}
\newcommand{\br}{\bar{\mathbb R}}
\newcommand{\Z}{{\mathbb Z}}

\newcommand{\R}{{\mathbb R}}
\newcommand{\E}{{\mathbb E}}

\renewcommand{\P}{{\mathbb P}}
\newcommand{\N}{{\mathbb N}}

\newcommand{\W}{{\cal W}}
\newcommand{\X}{{\cal X}}

\newcommand{\h}{{\cal H}}
\newcommand{\f}{{\cal F}}
\newcommand{\cd}{{\cal D}}

\renewcommand{\o}{\Pi}

\newcommand{\ttil}{\tilde t}
\newcommand{\totil}{\tilde t_0}
\newcommand{\atil}{\tilde a}

\newcommand{\etil}{\tilde \epsilon}

\newcommand{\util}{\tilde u}
\newcommand{\xtil}{\tilde x}
\newcommand{\xotil}{\tilde x_0}
\newcommand{\ytil}{\tilde y}

\newcommand{\pitil}{\tilde \pi}

\newcommand{\dbar}{\bar d}
\newcommand{\etahat}{\hat \eta}

\renewcommand{\a}{\alpha}
\renewcommand{\b}{\beta}
\newcommand{\g}{\gamma}
\newcommand{\G}{\Gamma}
\renewcommand{\k}{\kappa}

\renewcommand{\d}{\delta}

\newcommand{\e}{\epsilon}
\newcommand{\s}{\sigma}
\newcommand{\B}{{\cal B}}

\newcommand{\nn}{\nonumber}
\renewcommand{\=}{&=&}

\renewcommand{\le}{\leq}

\newcommand{\Xc}{\breve \X}
\newcommand{\muc}{\breve \mu}

\begin{document}

%\linespread{1.6}   %double spaced
%\normalsize        %set to normalsize so double spacing kicks in

\pagenumbering{roman}  % front matter number with roman

%%-------------------------------------Title page (unnumbered)
%\include{title}

\thispagestyle{empty}
\vspace*{0.5cm}
\begin{center}
\huge\bfseries
\thesistitle
\end{center}

\begin{center}
\vspace{2cm}

{\it \large \myname}\\

\vspace*{3.5cm}

\noindent
A dissertation submitted in partial fulfillment\\
of the requirements for the degree of \\
Doctor of Philosophy\\
Department of Mathematics\\
New York University\\
\thesismonth~\thesisyear\\
\end{center}

\vspace{1cm}

\begin{flushright}
   \makebox[2.2in]{\hrulefill}\\
   \vspace{-.2cm}
   \parbox[t]{1.8in}{\advisor}\\
\end{flushright}
\newpage

%%-------------------------------------Copyright  (unnumbered)
%%                                  (only if choose to copyright)
%\vspace*{2in}
%\begin{center}
%\copyright \hspace{.2cm} \myname\\
%All Rights Reserved, \thesisyear
%\end{center}

%\thispagestyle{empty} \hbox{   } \newpage

%%-------------------------------------BLANK (unnumbered)
%\thispagestyle{empty} \hbox{   } \newpage

%%-------------------------------------Frontispiece (optional, unnumbered)
%\include{frontispiece}

%%-------------------------------------Dedication, (first numbered)
\pagestyle{plain}
\vspace*{3in}
\begin{center}
{\large\em To my parents}
\end{center}

\addcontentsline{toc}{chapter}{Dedication}
\hbox{ } \newpage

%%-------------------------------------Acknowledgements
%\include{acknowledge}

\chapter*{Acknowledgements}
\addcontentsline{toc}{chapter}{Acknowledgements}

I am deeply indebted to my advisor Prof.~Charles M. Newman, who tolerated my earlier years of
idleness, patiently guided me through this work, and carefully read an earlier version of this
thesis and provided many corrections and suggestions. I especially thank Prof.~Krishnamurthi
Ravishankar for many helpful suggestions, who together with Chuck helped me overcome some of the
key difficulties. I would like to thank Prof.~S.~R.~S.~Varadhan, whose probability student seminar
benefited me greatly, and I thank both Prof.~Varadhan and Prof.~Gerard Ben Arous for reading a
draft of this thesis and pointing out some corrections. I would also like to thank my friends, Nikolaos
Zygouras, Toufic Suidan, Paolo Cascini, Fanfu Feng, Silas Alben and many others for many helpful and
entertaining conversations on both mathematics and other subjects, which made my graduate life much
more enjoyable and memorable. Lastly, I would like to thank my mother and brother for their support
during all the years I have been away and pursuing my studies, and my father who is no longer with
me today.

\newpage

%%-------------------------------------Preface
%\include{preface}
%\newpage

%%-------------------------------------Abstract
%\include{abstract}
\chapter*{Abstract}
\addcontentsline{toc}{chapter}{Abstract}

The Brownian Web (BW) is a family of coalescing Brownian motions starting
from every point in space and time $\R\times\R$. It was first introduced
by Arratia, and later analyzed in detail by T\'{o}th and Werner. More recently, Fontes,
Isopi, Newman and Ravishankar gave a characterization of the BW,
and general convergence criteria allowing either crossing or
noncrossing paths, which they verified for coalescing simple random walks.
Later Ferrari, Fontes, and Wu verified these criteria for a two dimensional
Poisson Tree. In both cases, the paths are noncrossing. In this thesis, we formulate
new convergence criteria for crossing paths, and verify them for non-simple coalescing
random walks satisfying a finite fifth moment condition. This is the first time
convergence to the BW has been proved for models with crossing paths. Several corollaries
are presented, including an analysis of the scaling limit of voter model interfaces that
extends a result of Cox and Durrett.

\newpage

%%-------------------------------------TOC
\tableofcontents
\newpage

%%-------------------------------------LOF
%\addcontentsline{toc}{chapter}{List of Figures}
%\listoffigures
%\newpage

%%-------------------------------------LOT
%\addcontentsline{toc}{chapter}{List of Tables}
%\listoftables
%\newpage

%%-------------------------------------List of appendices
%%                             (only use if more than one appendix)
%\listofappendices
%\newpage

%%-------------------------------------Thesis body

\pagenumbering{arabic}   %start arabic numbering

\chapter{Introduction and Results}
\setcounter{equation}{0}

The idea of the Brownian Web (BW) dates back to Arratia's thesis~\cite{kn:A0} in 1979, in
which he constructed a process of coalescing one-dimensional Brownian motions starting from every
point in space $\R$ at time zero. In a later unpublished manuscript~\cite{kn:A1}, Arratia
generalized this construction to a process of coalescing Brownian motions starting from
every point in space and time $\R\times\R$, which is essentially what is now often called
the Brownian Web. He also defined a dual family of backward coalescing Brownian motions
equally distributed (after a time reversal) with the BW which is now called the Dual Brownian Web.
Unfortunately, Arratia's manuscript was incomplete and never published, and the BW was not
studied again until a paper by T\'{o}th and Werner~\cite{kn:TW}, in which they constructed
and analyzed versions of the Brownian Web and Dual Brownian Web in great detail and used them
to construct a process they call the {\it True Self Repelling Motion}.

In both Arratia's and T\'{o}th and Werner's constructions of the BW, some semicontinuity
condition is imposed to guarantee a unique path starting from every space-time point. (For
example, in Arratia's~\cite{kn:A0} construction of coalescing Brownian motions starting from every
point on $\R$ at time 0, when multiple paths start from the same point, a unique path is chosen
so that if we regard the collection of paths as a mapping from $\R$ to the space of continuous
paths, then it is right continuous with left limits.) More
recently, Fontes, Isopi, Newman and Ravishankar~\cite{kn:FINR1, kn:FINR2} gave a different
formulation of the BW which provides a more natural setting for weak convergence, and they
coined the term Brownian Web. Instead of imposing semicontinuity conditions, multiple paths
are allowed to start from the same space-time point. Further, by choosing a suitable topology,
the BW can be characterized as a random variable taking values in a complete separable metric space
whose elements are compact sets of paths. In~\cite{kn:FINR2}, they gave general convergence
criteria allowing either crossing or noncrossing paths (i.e., almost surely, paths in the random
set may cross each other in the crossing paths case, and can never cross
each other in the noncrossing paths case), and they verified the criteria
for the noncrossing paths case for coalescing simple random walks. Recently, Ferrari, Fontes, and
Wu~\cite{kn:FFW} verified the same criteria for a two dimensional Poisson tree also with
noncrossing paths.

The main result of this thesis is the formulation of new convergence criteria for
the crossing paths case, and we verify these criteria for both discrete time and continuous
time one-dimensional coalescing {\em nonsimple} random walks satisfying a finite fifth moment
condition, which are models with crossing paths. The main technical distinction between the
noncrossing paths case and the crossing paths case is that if paths cannot cross, then they
form a totally ordered set, and one expects certain correlation inequalities to hold, which
is not the case if paths can cross.  We will also present some corollaries for one-dimensional
coalescing nonsimple random walks and their dual non-nearest-neighbor voter models.

In the remaining sections of this chapter, we will present some background on coalescing random
walks, their dual voter models, and the Brownian web. We will also state the main results of
this thesis.
\vspace{4mm}

\section{Coalescing random walks and voter models}

\noindent
Let $Y$, a random variable with distribution $\mu_Y$, denote the increment of an irreducible aperiodic
random walk on $\Z$. Throughout this thesis, all random walk increments will be distributed as $Y$.
We will always assume $\E[Y]=0$ and $\E[Y^2]=\s^2<+\infty$ unless a weaker hypotheses is explicitly
stated. For our main result, we will also need to assume $\E[|Y|^5] <+\infty$. For continuous time random
walks, we further assume that $\P(Y=0)=0$, and the random walk increments occur at rate 1.

\subsection{Coalescing random walks}

\noindent
The process of discrete time coalescing random walks with one walker starting from every space-time
lattice site is defined as follows. One walker starts at every site on the space-time lattice
$\Z\times\Z$ (first coordinate space, second coordinate time), and makes jumps at integer times
(including the time when it is born). All walkers have i.i.d. increments distributed as $Y$, and two
walkers move independently until they jump to the same site at the same time, then they coalesce.
The random walk paths are piecewise constant, right continuous with left limits, and have
discontinuities at integer times. A continuous version of the random walk path is obtained by linearly
interpolating between the consecutive space-time lattice sites from where the random walk makes its
jumps. Note that for non-simple random walks, two
interpolated random walk paths can cross each other many times before they actually
coalesce. If $Y$ was such that the random walks had period $d\neq 1$, as in the case of simple random
walks where $d=2$, then we would just have $d$ different copies of coalescing random walks on different
space-time sublattices, none of which interacts with the other copies. We will let $\X_1$ (with
distribution $\mu_1$) denote the random realization of such a collection of interpolated coalescing
random walk paths on $\Z\times\Z$, and let $\X_\d$ (with distribution $\mu_\d$, $0<\d<1$) be $\X_1$
rescaled with the usual diffusive scaling of $\d/\s$ in space, $\d^2$ in time. Sometimes we will
also need the piecewise constant version of $\X_\d$, which we denote by $\G_\d$, i.e., replace each
path in $\X_\d$ by its piecewise constant counterpart.

The continuous time analog of $\X_\d$, $\Xc_\d$ is defined as follows. One walker starts from every
point on $\Z \times \R$ and undergoes rate 1 jumps with increments distributed as $Y$.
The jump times are given by independent rate 1
Poisson clocks at each integer site $i\in\Z$. Two walkers coalesce when they first occupy the same site
at the same time. Clearly all walkers starting at the same site between two consecutive Poisson clock
rings will have coalesced by the time of the second clock ring. If we call the time and location at
which a Poisson clock rings a {\it jump point}, then the path of a random walk is piecewise constant
with discontinuites at the jump points. We can also define an interpolated version of the random walk
path, which consists of first a constant position line segment connecting the point of the walker's
birth to its first jump point, and then linearly interpolating between consecutive jump points along
its path. For random  walks born at a jump point, we will take two paths, one starting with a constant
position line segment, and one without. $\Xc_1$ (with distribution $\muc_1$) is then defined to be the
random variable consisting of all the interpolated coalescing random walk paths, and $\Xc_\d$ (with
distribution $\muc_\d$) is $\Xc_1$ diffusively rescaled. We will denote the piecewise constant
version of $\Xc_\d$ by $\breve\G_\d$.

For a system of d-dimensional ($d\geq 1$) discrete time coalescing random walks starting from a
space-time subset
$A \subset \Z^d\times\Z$, we will denote the set of occupied sites in $\Z^d$ at time $n$ by $\xi^A_n$.
(To be consistent with our definition of $\X_1$, we also define the random walks so that they jump
at the time of birth. But the random walks' positions at time $n$ is now taken to be the positions of
the interpolated random walk paths at time $n$.) In the special case when $A = B\times\{n_0\}$ for
some $B\subset\Z^d$, we will denote it by $\xi^{B,n_0}_n$; and when $n_0=0$, we simply denote it by
$\xi^B_n$. For most of this thesis, we will only deal with one-dimensional coalescing random
walks. To simplify the notation, we will use $\xi^A_s$ also to denote the continuous time analogue
of $\xi^A_n$, where $A\subset\Z^d\times\R$ and the random walks jump with rate 1 and increments
distributed as $Y$.

The main result of this thesis (see Theorem \ref{teo:walkconv} below) is that if
$\E[|Y|^5]<+\infty$, then $\X_\d$ and $\Xc_\d$ converge in distribution to the BW as $\d\to 0$.

\subsection{Voter models}\label{voter models}

\noindent
The voter model was first introduced in the probability literature by Clifford and
Sudbury~\cite{kn:CS}, and Holley and Liggett~\cite{kn:HL}. In population genetics, a variant
of the voter model, the stepping stone model, was first introduced by Kimura~\cite{kn:Kimu}.
A two type d-dimensional discrete time voter model is defined as follows. The state space is
$\{0,1\}^{\Z^d}$ with product topology. A metric that generates the product topology is
$\|\eta-\zeta\| = \sum_{x\in\Z^d} 2^{-|x|_{\infty}}|\eta(x)-\zeta(x)|$. We will denote the
state at time $n$ of the voter model by $\phi_n$ (with distribution $\nu_n$). $\phi_n(x)$, the value
of $\phi_n$ at site $x$, can be regarded as the opinion, or political affiliation,
of the voter at site $x$ at time $n$, hence the name {\em voter model}. The initial configuration
is specified by $\phi_0$. If we also use $Y$ to denote the increment of a d-dimensional random
walk, then at each time $n\geq 1$, we update $\phi_n(x)$ by setting $\phi_n(x) = \phi_{n-1}(x+Y_{x,n})$,
where $\{Y_{x,n}\}_{x\in\Z^d, n\in\N}$ are i.i.d. $\Z^d$-valued random variables distributed as $Y$.
Given any realization of $\{Y_{x,n}\}_{x\in\Z^d, n\in\N}$, if we trace backward in time the geneaologies
of the opinions of all sites at all times, then the collection of all such geneaology trees is exactly
a realization of coalescing random walk paths running backward in time, with one walker starting at every
site in $\Z^d$ at every positive integer time. This provides a natural coupling between coalescing
random walks and voter models, and leads to the duality relation
\beq \label{duality}
\P[\phi_n(x)=1 \hbox{ for some }(x,n) \in A] = \P[\phi_0(y)=1 \hbox{ for some } y \in \hat\xi^A_0],
\eeq
where $A\subset \Z^d\times\Z^+$ is a set of space-time lattice sites with positive time, and
$\hat\xi^A_0$ is the set of occupied sites at time $0$ of a system of coalescing random walks running
backward in time, with one walker starting at every site in $A$. For more details on the voter
model, see, e.g., \cite{kn:HL,kn:L}.

The definition above easily generalizes to a multitype voter model with state space  $S^{\Z^d}$,
where $S$ can be any finite or countable set. For continuous time, the opinion
at each site is updated independently according to a rate 1 Poisson clock. Whenever a clock rings
at a site $x$, another site $y\in\Z^d$ is picked with $(y-x)$ distributed as the random walk increment
$Y$, and the opinion at site $x$ changes to that of site $y$. The dual of this continuous time voter
model is the process of rate 1 continuous time coalescing random walks running backward in time,
with one walker starting at every site in $\Z^d\times\R^+$. The two models are coupled through the
realization of the jump points and the random walk increment associated with each jump point.
We will denote the continuous time voter model also by $\phi_t$ (with distribution $\nu_t$).
The duality relation $(\ref{duality})$ also holds for continuous time.

\section{Brownian web}
\setcounter{equation}{0}

\noindent
One way of thinking about the Brownian web is to view it as the diffusive scaling limit of the system
of coalescing random walks $\X_1$. Heuristically, we expect to obtain as the limit a set of coalescing
Brownian motions with one Brownian motion starting from every space-time point. The main conceptual
difficulty is that there are uncountable number of space-time points, and in general we would like to
deal only with a countable number of Brownian motions because of the countable additivity of measures.
Fortunately, it turns out that the limiting object (the Brownian web) is fully determined by a countable
number of coalescing Brownian motions starting from a countable dense subset $\cal D$ of $\R^2$ (see
Theorem \ref{teo:char} below). Coalescing Brownian motion paths starting from space-time points
off the dense set $\cal D$ can be constructed by taking limits of paths starting from $\cal D$ using
the noncrossing property of coalescing Brownian motions. See Arratia~\cite{kn:A1} and T\'{o}th and
Werner~\cite{kn:TW} for two such constructions. The Brownian web intrinsically contains space-time
points from which multiple paths start out. In both Arratia's and T\'{o}th and Werner's construction,
only one path
is retained starting from every space-time point satisfying some semicontinuity conditions. Here we
follow a different approach by Fontes, Isopi, Newman and Ravishankar~\cite{kn:FINR1,kn:FINR2}, who
coined the term Brownian Web. Their approach is to regard the Brownian web as a random variable taking
values in the space of sets of paths. Thus there is no need to throw away paths when multiplicity
arises. Another advantage is that, by choosing
the topology approriately for the space of sets of paths, the Brownian web takes values in a complete
separable metric space, thus providing a natural setting for establishing weak convergence results,
which is the central theme of this thesis.

We now recall Fontes, Isopi, Newman and Ravishankar's \cite{kn:FINR1,kn:FINR2} choice of the metric
space in which the Brownian Web takes its values.

Let $(\bar\R^2,\rho)$ be the completion (or compactification) of $\R^2$
under the metric $\rho$, where
\begin{equation}
\label{rho}
\rho((x_1,t_1),(x_2,t_2))=
\left|\frac{\tanh(x_1)}{1+|t_1|}-\frac{\tanh(x_2)}{1+|t_2|}\right|
\vee|\tanh(t_1)-\tanh(t_2)|.
\end{equation}
$\bar\R^2$ can be thought of as the image of $[-\infty,\infty]\times[-\infty,\infty]$
under the mapping
\begin{equation}
\label{compactify}
(x,t)\leadsto(\Phi(x,t),\Psi(t))
\equiv\left(\frac{\tanh(x)}{1+|t|},\tanh(t)\right).
\end{equation}
We can think of the mapping as first squeeze $\bar\R^2$ to the square $[-1,1]\times[-1,1]$ by the
mapping $(\tanh x, \tanh t)$, and then the $x$ coordinate is squeezed even further depending on its
time coordinate such that the top and bottom edge of the square is squeezed to two points $(0,\pm 1)$.

For $t_0\in[-\infty,\infty]$, let $C[t_0]$ denote the set of functions
$f$ from $[t_0,\infty]$ to $[-\infty,\infty]$ such that $\Phi(f(t),t)$
is continuous. Then define
\begin{equation}
%\label{omega}
\o=\bigcup_{t_0\in[-\infty,\infty]}C[t_0]\times\{t_0\},
\end{equation}
where $(f,t_0)\in\o$ represents a path in $\br^2$ starting at
$(f(t_0),t_0)$.
For$(f,t_0)$ in $\o$, we denote by $\hat f$ the function that extends $f$
to all
$[-\infty,\infty]$ by setting it equal to $f(t_0)$ for $t<t_0$. Then we
take
\begin{equation}
\label{d}
d((f_1,t_1),(f_2,t_2))=
(\sup_t|\Phi(\hat{f_1}(t),t)-\Phi(\hat{f_2}(t),t)|)
\vee|\Psi(t_1)-\Psi(t_2)|,
\end{equation}
i.e., after applying the mapping $(\Phi, \Psi)$ to $(\hat{f_1}(t),t)$ and $(\hat{f_2}(t),t)$, the
distance between the two original paths is then taken to be the maximum of (i) the supremum norm
distance between the two image paths, and (2) the absolute difference in the starting times of the two
image paths. Note that $(\o,d)$ is a complete separable metric space.

Let now $\h$ denote the set of compact
subsets of $(\o,d)$, with $d_\h$ the induced Hausdorff metric, i.e.,
\begin{equation}
\label{dh}
d_\h(K_1,K_2)=\sup_{g_1\in K_1}\inf_{g_2\in K_2}d(g_1,g_2)\vee
                 \sup_{g_2\in K_2}\inf_{g_1\in K_1}d(g_1,g_2).
\end{equation}
$(\h,d_\h)$ is also a complete separable metric space. Let $\f_\h$
denote the Borel $\s$-algebra generated by $d_\h$.

In~\cite{kn:FINR1,kn:FINR2}, the Brownian Web ($\bar\W$ with measure $\mu_{\bar\W}$) is constructed
as a $(\h, \f_\h)$ valued random variable, with the following characterization theorem.
\bteo
\label{teo:char}
There is an \( ({\cal H},{\cal F}_{{\cal H}}) \)-valued random variable \(
\bar{\W} \)
whose distribution is uniquely determined by the following three
properties.
\begin{itemize}
         \item[(o)]  from any deterministic point \( (x,t) \) in
$\r^{2}$,
         there is almost surely a unique path \( {W}_{x,t} \)
starting
         from \( (x,t) \).

         \item[(i)]  for any deterministic \( n, (x_{1}, t_{1}), \ldots,
         (x_{n}, t_{n}) \), the joint distribution of \(
         {W}_{x_{1},t_{1}}\),\!\!\! \(\ldots, {W}_{x_{n},t_{n}} \) is that
         of coalescing Brownian motions (with unit diffusion constant),
and

         \item[(ii)]  for any deterministic, dense countable subset  \(
{\cal
         D} \) of \( \r^{2} \), almost surely, \( \bar{\W} \) is the
closure in
         \( ({\cal H}, d_{{\cal H}}) \) of \( \{ {W}_{x,t}: (x,t)\in
         {\cal D} \}. \)
\end{itemize}
\eteo
The $(\h, \f_\h)$-valued random variable $\bar\W$ in Theorem $\ref{teo:char}$ is
called the standard Brownian Web.

The Brownian web $\bar\W$ uniquely determines a dual (backward) Brownian web $\hat{\bar\W}$, which
is equally distributed with the standard Brownian web $\bar\W$ except for a time reversal.
The pair $(\bar\W, \hat{\bar\W})$ forms what is called the {\em double Brownian web} $\bar\W^D$
with the property that, almost surely, paths in $\bar\W$ and $\hat{\bar\W}$ {\em reflect} off each
other and never cross. The double Brownian web is most easily recognized as the limit of coalescing
simple random walks. For coalescing simple random walks, only walks starting from lattice sites
$(x,m)\in\Z^2$ with $x+m$ having the same parity (even or odd) interact with each other. Take the set
of walks that start from $(x,m)$ with $x+m$ even. Then any realization of such forward coalescing
random walk paths uniquely determines a set of backward (in time) coalescing simple random walk paths
with one path starting from every site $(y,n)\in\Z^2$ with $y+n$ odd, provided that we require the
backward paths never cross the forward paths. The resulting system of backward coalescing simple
random walks is equally distributed with the forward coalescing system except for a time reversal.
Under the diffusive scaling limit, it is then seen that the joint distribution of the forward and
backward systems of coalescing simple random walks converge in distribution to the double Brownian
web. For more on the Brownian web and the double Brownian web, see~\cite{kn:A0, kn:A1, kn:TW, kn:STW,
kn:FINR1, kn:FINR2}.

\section{Convergence criteria and main result}
\setcounter{equation}{0}

\noindent
In~\cite{kn:FINR2}, a set of general convergence criteria were formulated for measures supported
on compact sets of paths which can cross each other. However, one of the conditions $(B_2')$
turns out to be difficult to verify for the coalescing nonsimple random walks $\{\X_\d\}$ and
$\{\Xc_\d\}$.
In our modified general convergence theorem \ref{teo:mgconv}, we will replace condition $(B_2')$
by an alternative condition $(E_1)$, and we will verify $(E_1)$ and the other convergence criteria for
$\{\X_\d\}$ and $\{\Xc_\d\}$ under the assumption that $\E[|Y|^5]<+\infty$.

We first recall the convergence criteria formulated in~\cite{kn:FINR2} (for the crossing
paths case) for a family of $(\h,\f_\h)$-valued random variables $\{X_n\}$ with distributions
$\{\mu_n\}$.
\vspace{2mm}

\noindent $(I_1)$ There exist single path valued random variables
$\theta^y_n\in X_n$, for $y\in\R^2$, satisfying:
for $\cd$ a deterministic countable dense subset of $\R^2$, for any deterministic
$y^1,\ldots,y^m\in\cd$, ${\theta^{y^1}_n,\ldots,\theta^{y^m}_n}$
converge jointly in distribution as $n\to +\infty$ to coalescing Brownian
motions (with unit diffusion constant) starting at $y^1,\ldots,y^m$.

Let $\Lambda_{L,T} = [-L,L]\times [-T,T] \subset \R^2$.
For $x_0,t_0\in\r$ and $u,t>0$, let $R(x_0,t_0;u,t)$ denote the
rectangle $[x_0-u,x_0+u]\times[t_0,t_0+t]$ in $\R^2$. Define
$A_{t,u}(x_0,t_0)$ to be the event (in $\f_\h$) that $K$ (in $\h$)
contains a path touching both $R(x_0,t_0;u,t)$ and (at a
later time) the left or right boundary of the bigger rectangle
$R(x_0,t_0;17u,2 t)$ (the number 17 is chosen to avoid fractions later).
Then the following is a tightness condition for $\{X_n\}$: for every
$u, L, T\in (0, +\infty)$,
\beqnn
\!\!\!\!\!\!\!\!\!\!\!
(T_1)&&\!\!\!\!\!\! \tilde g(t,u;L,T)\equiv t^{-1}\limsup_{n\to +\infty}\,\,
                \sup_{(x_0,t_0) \in \Lambda_{L,T}}
            \mu_n(A_{t,u}(x_0,t_0))\to0\mbox{ as }t\to 0^+ \ .
\eeqnn
As shown in \cite{kn:FINR2}, if $(T_1)$ is satisfied, one can construct compact
sets $G_\e \subset \h$ for each $\e>0$, such that $\mu_n(G_\e^c) < \e$ uniformly in $n$.
$G_\e$ consists of compact subsets of $\Pi$ whose image under the map $(\Phi, \Psi)$
are equicontinuous with a modulus of continuity dependent on $\e$.

For $K\in \h$ a compact set of paths in $\Pi$, define the counting variable
${\cal N}_{t_0,t}([a,b])$ for $a,b,t_0,t\in \R, a<b, t>0$ by
\begin{eqnarray}
{\cal N}_{t_0,t} ([a,b])
\= \{ y \in \R\,|\,\exists\, x \in [a,b]
\hbox{ and a path in } K \hbox{ which touches} \nn \\
&&\hspace{1.6cm}\hbox{both }  (x,t_0) \hbox{ and }
(y,t_0 + t)\}. \label{calN}
\end{eqnarray}
Let $l_{t_0}$ (resp., $r_{t_0}$) denote the leftmost (resp., rightmost) value in $[a,b]$
with some path in $K$ touching $(l_{t_0}, t_0)$ (resp., $(r_{t_0}, t_0)$). Also
define ${\cal N}^+_{t_0,t} ([a,b])$ (resp. ${\cal N}^-_{t_0,t} ([a,b])$) to be
the subset of ${\cal N}_{t_0,t} ([a,b])$ due to paths in $K$ that touch $(l_{t_0}, t_0)$
(resp., $(r_{t_0}, t_0)$). The last two conditions for the convergence of
$\{X_n\}$ to the Brownian Web are

\beqnn
\!\!\!\!
(B'_1)&&\!\!\!\!\!\!\!\!  \forall\b >0,
\limsup_{n\to +\infty}\sup_{t >\b}  \sup_{t_0,a \in \R} \mu_n
(|{\cal N}_{t_0,t} ([a-\e,a+\e]) | > 1) \to 0 \hbox{ as } \e
\to 0^+,
\eeqnn
\beqnn
(B'_2)&&\!\!\!\!\!\!\!\!  \forall\b >0,
\frac{1}{\e} \limsup_{n\to +\infty}\sup_{t >\b}
\sup_{t_0,a \in \R} \mu_n (
{\cal N}_{t_0,t} ([a-\e,a+\e])\ne {\cal N}^+_{t_0,t}([a-\e,a+\e])\\
&&\hspace{3.8cm}
\cup{\cal N}^-_{t_0,t}([a-\e,a+\e]))\to0\hbox{ as }\e\to0^+.
\eeqnn

The general convergence theorem of~\cite{kn:FINR2} is the following,
\bteo \label{teo:gconv}
Let $\{X_n\}$ be a family of $(\h, \f_\h)$ valued random variables satisfying
conditions $(I_1), (T_1), (B'_1)$ and $(B'_2)$, then $X_n$ converges in distribution
to the standard Brownian Web $\bar\W$.
\eteo

Condition $(B'_1)$ guarantees that for any subsequential limit $X$ of $\{X_n\}$ (with distribution
$\mu_X$), and for any deterministic point $y\in \R^2$, there is $\mu_X$ almost surely at most one
path starting from $y$. Together with condition $(I_1)$, this implies that for a
deterministic countable dense set $\cd \subset \R^2$, the distribution of paths in $X$
starting from finite subsets of $\cd$ is that of coalescing Brownian motions. This shows $X$
contains at least as many paths as the Brownian web $\bar\W$. Conditions
$(B'_1)$ and $(B'_2)$ together imply that for the family of counting random variables
$\eta(t_0,t;a,b) = |{\cal N}_{t_0,t}([a,b])|$, we have
$\E[\eta_X(t_0,t;a,b)] \leq \E[\eta_{\bar\W}(t_0,t;a,b)] = 1 +\fr{b-a}{\sqrt{\pi t}}$
for all $t_0,t,a,b\in\R$ with $t>0, a<b$. By Theorem $4.6$ in~\cite{kn:FINR2} and the remark following
it, this fact implies that $X$ contains no extra paths besides the Brownian web $\bar\W$, thus
$X$ is equidistributed with $\bar\W$. For the systems of coalescing random walks $\{\X_\d\}$ and
$\{\Xc_\d\}$,
we have not yet been able to verify condition $(B'_2)$, but an examination of the proof of Theorem
$4.6$ in~\cite{kn:FINR2}  shows that we can also use the dual family of counting random
variables
\begin{eqnarray}
\etahat_X(t_0,t;a,b) = |\{x\in(a,b) &|& \exists \hbox{ a path in } X \hbox{ touching}
\label{etahat} \\
&&\hbox{both } \R\times\{t_0\} \hbox{ and } (x, t_0+t)\} |.  \nn
\end{eqnarray}
By a duality argument~\cite{kn:TW} (see also~\cite{kn:A0, kn:A1, kn:FINR2}), $\etahat$
and $\eta-1$ are equally distributed for the Brownian Web $\bar\W$. We can then replace
$(B'_2)$ by
\vskip 0.3cm
\noindent$(E_1)$ If $X$ is any subsequential limit of $\{X_n\}$, then
$\forall t_0,t,a,b\in\R$ with $t>0$ and $a<b$,
$\E[\etahat_X(t_0,t;a,b)] \leq \E[\etahat_{\bar\W}(t_0,t;a,b)] = \fr{b-a}{\sqrt{\pi t}}$.
\vspace{2mm}

\noindent With this change, we obtain our modified general convergence theorem,
\bteo \label{teo:mgconv}
Let $\{X_n\}$ be a family of $(\h, \f_\h)$ valued random variables satisfying
conditions $(I_1), (T_1), (B'_1)$ and $(E_1)$, then $X_n$ converges in distribution
to the standard Brownian Web $\bar\W$.
\eteo
The main result of this paper is
\bteo \label{teo:walkconv}
If the random walk increment $Y$ satisfies $\E[|Y|^5]< +\infty$, then
$\{\X_\d\}$ and $\{\Xc_\d\}$ satisfy the conditions of Theorem $\ref{teo:mgconv}$, and hence
converge in distribution to $\bar\W$.
\eteo
\brm \label{rmk:1.6}
Condition $(E_1)$ in our general convergence theorem \ref{teo:mgconv} may seem
strong, but as we will show in our proof of $(E_1)$ for $\{\X_\d\}$ and $\{\Xc_\d\}$ in Section
\ref{E1}, the key ingredients are the Markov property of the random walks and an upper bound of
the type  $\limsup_{\d\downarrow 0}\E[\etahat_{\X_\d}(t_0, t;a,b)] \leq C$ for some finite constant
$C$ depending on $t,a,b$.
\erm

\chapter{Random Walk Estimates} \label{randomwalk}
\setcounter{equation}{0}
In this chapter, we first introduce some notation, and then list some basic facts about
random walks that will be used throughout the rest of the thesis. Once acquainted with the basic
notation, the reader may skip the rest of the chapter until the results in this chapter are referred to.

Given macroscopic space and time coordinates $(x,t)\in\R^2$, define their microscopic counterparts
before diffusive scaling by $\ttil = t\d^{-2}$ and $\xtil = x\s\d^{-1}$. Quantities
such as $\util, \totil$ are defined from $u, t_0$ similarly depending on whether they are space or
time units. Since $\mu_\d$ (resp., $\muc_\d$) and
$\mu_1$ (resp., $\muc_1$) are related by diffusive scaling, we will do most of our analysis using
$\mu_1$ (resp., $\muc_1$), with $x,t, u,t_0$ for $\mu_\d$ (resp., $\muc_\d$) replaced by
$\xtil,\ttil,\util,\totil$ for $\mu_1$ (resp., $\muc_1$).

For both discrete time and continuous time, we will denote the piecewise constant version of
the path of a random walk (which by definition is right continuous with left limits in time)
starting at position $x$ at time $t_0$ by $\pi^{x,t_0}(s)$. We will denote the linearly interpolated
version of the random walk path by $\k^{x,t_0}(s)$. Denote the event that the path of a random walk
$\pi^{x,t_0}(s)$ (either discrete or continuous time) starting at $(x,t_0)$ stays inside the interval
$[a,b]$ containing $x$ up to time $t$ by $B^{x,t_0}_{[a,b],t}$.

Given any $r\in\Z$, we also define the following stopping times associated with either a discrete or
a continuous time random walk $\pi^{x,t_0}$
\begin{eqnarray}
&& \tau_r^{x,t_0} = \inf \{t\geq t_0 \ |\ \pi^{x,t_0}(t) = r\}, \label{tau} \\
&& \tau_{r^+}^{x, t_0} = \inf \{t\geq t_0 \ |\ \pi^{x,t_0}(t) \geq r\}. \label{tau+}
\end{eqnarray}
When the time coordinate in the superscripts of $\pi, \k, B, \tau, \tau_{^+}$ is 0, we will suppress
it. We will use $\P_x$ and $\E_x$ to denote probability and expectation for a random walk process
starting from $x$ at time 0. $\P_{x,y}$ and $\E_{x,y}$ will correspond to two independent random walks
starting at $x$ and $y$ at time 0.

The following lemmas are stated for random walk paths $\pi^x$ and $\pi^y$, which can be interpreted
for both discrete and continuous time. The random walk increment $Y$ is always assumed to be that of
an irreducible and aperiodic random walk with $\E(Y)=0$ and $\E[Y^2]<+\infty$, unless a different
moment condition is explicitly stated. For continuous time random walks, we always assume it jumps
with rate 1 unless otherwise explicitly stated.

\blem \label{lem:excursion}
Let $\pi^x, \pi^y$ be two independent random walks with increment $Y$ starting at $x, y\in\Z$ at time 0.
Let $\tau_{x,y} = \inf \{ t\geq 0\, |\, \pi^x(t) = \pi^y(t) \}$, which is a stopping time, and let
$l(x,y) = \sup_{t\in[0,\tau_{x,y}]}|\pi^x(t)-\pi^y(t)|$. Then $\tau_{x,y}$ and $l(x,y)$
are almost surely finite.
\elem
{\bf Proof.} Let $\bar\pi^{y-x}(t) = \pi^y(t) - \pi^x(t)$. Then $\bar\pi^{y-x}$ is an irreducible
aperiodic symmetric random walk starting at $y-x$. For discrete time, $\bar\pi^{y-x}$ has increment
distributed as $\mu_Y\ast\mu_{-Y}$; for continuous time,  $\bar\pi^{y-x}$ is a rate 2 random walk
with increment distributed as $\fr{1}{2}\mu_Y +\fr{1}{2}\mu_{-Y}$. The lemma is simply a consequence
of the recurrence of $\bar\pi^{y-x}$, which requires less than finite second moment of $Y$.

\blem \label{lem:hitprob}
Let $\pi^x, \pi^y$, $\tau_{x,y}$ be as in Lemma $\ref{lem:excursion}$. Then
$\P_{x,y}(\tau_{x,y} > t) \leq \fr{C}{\sqrt{t}}|x-y|$ for some constant $C$ independent of $t,x$
and $y$.
\elem
{\bf Proof.} Let $\bar\pi^{y-x}(t) = \pi^y(t) - \pi^x(t)$ as in the proof of
Lemma $\ref{lem:excursion}$. Let $\bar \P_{y-x}$ denote probability for this random
walk, and let $\bar\tau^{y-x}_0$ denote the stopping time when the random walk $\bar\pi^{y-x}$ first
lands at the site 0. Then $\P_{x,y}(\tau_{x,y} >t) = \bar\P_{y-x}(\bar\tau^{y-x}_0 >t)$. When
$|x-y|=1$, it is a standard fact (see, e.g., Proposition 32.4 in~\cite{kn:S})
that this probability is bounded by $\fr{C}{\sqrt{t}}$. When $|x-y|>1$, we can without loss of
generality assume $x<y$ and regard $\{\pi^x, \pi^y\}$ as a subset of the system of
coalescing random walks $\xi^{\{x, x+1, ...,y\}}$ up to time $\tau_{x,y}$. Then
\beqnn
&&\P_{x,y}(\tau_{x,y} >t) \leq \P(|\xi^{\{x,...,y\}}_t|>1) \\
&=& \P(\bigcup_{i=x}^{y-1}\{\tau_{i,i+1}>t\})
\leq (y-x)\P_{0,1}(\tau_{0,1}>t)\leq \fr{C(y-x)}{\sqrt{t}},
\eeqnn
which establishes the lemma.

\blem \label{lem:bar}
Let $u>0$ and $t>0$ be fixed, and let $\pi(s) = \pi^{0,0}(s)$
be a random walk starting from the origin at time 0. Let $\util, \ttil$ and
the event $B^{0,0}_{[-\util,\util], \ttil}$ be defined as at the beginning of this
chapter (note that they depend on $\d$), and let $(B^{0,0}_{[-\util,\util], \ttil})^c$ be the
complement of $B^{0,0}_{[-\util,\util], \ttil}$. If ${\cal B}_s$ is a standard Brownian motion
starting from 0, then
\beqnn
0<\lim_{\d\to 0^+} \P_0[(B^{0,0}_{[-\util,\util], \ttil})^c] =
 \P(\sup_{s\in[0,t]}|{\cal B}_s| > u)
                              < 4e^{-\fr{u^2}{2t}}. \nn
\eeqnn
\elem
{\bf Proof.}
The limit follows from Donsker's invariance principle for random walks. The
first inequality is trivial, and the second inequality follows from a well-known
computation for Brownian motion using the reflection principle.

\blem \label{lem:crossbar}
Let $u, t, \util, \ttil$ be as before. Let $\pi^x, \pi^y$ and $\tau_{x,y}$ be as in Lemmas
$\ref{lem:excursion}-\ref{lem:hitprob}$ with $x<y$. Let
$\tau_{x,y,\util^+} = \inf \{t\geq 0\,|\, \pi^x(t) -\pi^y(t)\geq\util\}$.
If $\E[|Y|^3]<+\infty$, and $\d$ is sufficiently small, then we have
\beqnn
    \P_{x,y}(\tau_{x,y,\util^+}<(\tau_{x,y} \wedge \ttil)) < C(t,u)\d
\eeqnn
for some constant $C(t,u)$ depending only on $t$ and $u$.
\elem
{\bf Proof.} Let $z=x-y<0$. Note that $x,y,z$ are fixed while $\util,\ttil\to+\infty$
as $\d\to 0$. For the difference of the two walks
$\bar\pi^z(s) = \pi^x(s) - \pi^y(s)$, we denote the first time when $\bar\pi^z(s)=0$ by
$\bar\tau_0^z$, and the first time when $\bar\pi^z(s) \geq \util$ by $\bar\tau_{\util^+}^z$. We are
using the bar $\bar\cdot$ to emphasize the fact that we are studying the symmetrized random walks.
The inequality then becomes
\begin{eqnarray}
\bar\P_z(\bar\tau_{\util^+}^z <(\bar\tau_0^z\wedge\ttil)) < C(t,u)\d.  \label{ineq1}
\end{eqnarray}
For simplicity, we will only prove $(\ref{ineq1})$ for the discrete time case. For the continuous
time case, only the notations will be different. We will first prove that, for $\d$ sufficiently small,
\begin{eqnarray}
\bar\P_w(\bar\tau_{\util^+}^w <(\bar\tau_0^w\wedge\ttil)) < C'(t,u)|w|\d. \label{ineq2}
\end{eqnarray}
By the strong Markov property,
\begin{eqnarray}
&&      \bar\P_w (\bar\tau_0^w >\ttil) \nn \\
&\geq&  \sum_{k=\lceil\util\rceil}^{+\infty}
        \bar\P_w(\bar\tau_{\util^+}^w<(\bar\tau_0^w\wedge\ttil),
        \bar\pi^w(\bar\tau_{\util^+}^w)=k,
        B^{k, \bar\tau^w_{\util^+}}_{[k-\fr{\util}{2},k+\fr{\util}{2}],\ttil}) \nn\\
&=&     \sum_{k=\lceil\util\rceil}^{+\infty}
        \sum_{n=0}^{\lfloor\ttil\rfloor}
        \bar\P_w(\bar\tau_{\util^+}^w=n, n<\bar\tau_0^w,\bar\pi^w(n)=k)
                      \bar\P_k(B^k_{[k-\fr{\util}{2}, k+\fr{\util}{2}],\ttil-n})
    \nn\\
&\geq&  \sum_{k=\lceil\util\rceil}^{+\infty}
        \sum_{n=0}^{\lfloor\ttil\rfloor}
        \bar\P_w(\bar\tau_{\util^+}^w=n, n<\bar\tau_0^w,\bar\pi^w(n)=k)
                      \bar\P_k(B^k_{[k-\fr{\util}{2}, k+\fr{\util}{2}],\ttil})
    \nn\\
&\geq&  C''(t,u)\bar\P_w(\bar\tau^w_{\util^+} <(\bar\tau_0^w\wedge\ttil)). \nn
\end{eqnarray}
If $\d$ is sufficiently small, the last inequality is valid by Lemma $\ref{lem:bar}$.
Also by Lemma $\ref{lem:hitprob}$,
\beqnn
\bar\P_w (\bar\tau^w_0 >\ttil) < \fr{C}{\sqrt{\ttil}}|w| = \fr{C|w|}{\sqrt{t}}\d,
\eeqnn
together they give $(\ref{ineq2})$.

To show $(\ref{ineq1})$, we condition at the first time when $\bar\pi^z(s)\geq 0$, which we denote
by $\bar\tau^z_{0^+}$. Then by the strong Markov property and $(\ref{ineq2})$,
\begin{eqnarray}
&&      \bar\P_z(\bar\tau^z_{\util^+}<(\bar\tau^z_0\wedge\ttil)) \nn\\
&=&     \sum_{w=1}^{+\infty}
        \sum_{n=0}^{\lfloor\ttil\rfloor}
        \bar\P_z[\bar\tau^z_{0^+}=n, \bar\pi^z(n)=w]
        \bar\P_w[\bar\tau^w_{\util^+} <(\bar\tau^w_0\wedge(\ttil-n))] \nn\\
&<&     \sum_{w=1}^{+\infty}
        \sum_{n=0}^{\lfloor\ttil\rfloor}
        \bar\P_z[\bar\tau^z_{0^+}=n, \bar\pi^z(n)=w]\ C'(t,u)|w|\ \d
    \nn\\
&<&     C'(t,u)\ \d\ \bar\E_z[\bar\pi^z(\bar\tau^z_{0^+})] < C(t,u) \d.  \nn
\end{eqnarray}
The last inequality follows from our assumption $\E[|Y|^3]<+\infty$ and the following
two lemmas.

\blem \label{lem:overshootlimit}
Let $\pi^x$ be a random walk with increment $Y$ starting from $x<0$ at time 0.
If $\E[Y^2] < +\infty$, then the overshoot $\pi^x(\tau_{0^+}^x)$ has a limiting distribution
as $x\to -\infty$. In terms of the ladder variable $Z = \pi^{0}(\tau^0_{1^+})$,
\beqnn
\lim_{x\to -\infty} \P[\pi^x(\tau^x_{0^+})=k] = \fr{\P[Z\geq k+1]}
{\E[Z]}.
\eeqnn
\elem
{\bf Proof.}This is a standard fact from renewal theory, see e.g. Proposition 24.7
in~\cite{kn:S}.

\blem\label{lem:overshoot}
Let $\pi^x$, $Y$ and $Z$ be as in the previous lemma. If $\E[|Y|^{r+2}]<+\infty$ for some
$r>0$, then $[\pi^x(\tau_{0^+}^x)]^r$ is uniformly integrable in $x\in\Z_-$, and
\beqnn
\lim_{x\to -\infty} \E\big[ [\pi^x(\tau_{0^+}^x)]^r\big]
    = \fr{1}{\E[Z]}\sum_{k=1}^{+\infty}k^r\P[Z\geq k+1] <+\infty.
\eeqnn
\elem
{\bf Proof.} We may assume $\pi^x$ is a discrete time random walk, since the continuous time
random walk is just a random time change of the discrete time walk, which does not change
the overshoot distribution. Note that if we let $\g^x$ denote the discrete time random walk
starting from $x<0$ at time 0 with increment distributed as $Z$, then $\g^x$ simply records
the successive maxima of the random walk $\pi^x$, so the overshoots $\pi^x(\tau_{0^+}^x)$
and $\g^x(\tau_{0^+}^x)$ are equally distributed. By a last passage decomposition for $\g^x$,
\beqnn
\P[\g^x(\tau_{0^+}^x)=k] = \sum_{i=x}^{-1}G_\g(x, i)\P[Z=k-i] \leq \P[Z\geq k+1],
\eeqnn
where $G_\g(x,i)$ is the probability $\g^x$ will ever visit $i$. Since
$\E[|Y|^{r+2}]<+\infty$ implies $\E[Z^{r+1}]<+\infty$ (see e.g. problem 6 in Chapter
IV of~\cite{kn:S}), we have
$\E\big[ [\pi^x(\tau_{0^+}^x)]^r\big]\leq\sum_{k=1}^{+\infty}k^r\P[Z\geq k+1] <+\infty$,
giving uniform integrability. The rest then follows from Lemma $\ref{lem:overshootlimit}$ and
dominated convergence.

\blem\label{lem:densityupperbound}
Let $\xi^{\Z}_t$ be a system of coalescing random walks (either discrete or continuous time)
starting from every site of $\Z$ at time 0, whose random walk increments are distributed as $Y$
with $\E[Y^2]<+\infty$. Then
$p_t\equiv \P(0\in \xi^{\Z}_t) \leq \fr{C}{\sqrt t}$ for some constant $C$ independent of
the time $t$.
\elem
\brm
We present two proofs, the first of which works for both discrete and continuous time, and
is an adaptation of the argument used by Bramson and Griffeath~\cite{kn:BG} to establish similar
upper bounds for continuous time coalescing simple random walks in $Z^d, d\geq 2$. The second
proof is special to continuous time walks, and can also be found in the paper of Bramson and
Griffeath~\cite{kn:BG}. In Corollary \ref{cor:density} below, we will prove that in fact
$p_t \sim 1/(\sigma \sqrt{\pi t})$ as $t\to +\infty$ under the stronger assumption that
$\E[|Y|^3]<+\infty$.
\erm
{\bf First Proof.} Let $B_M = [0, M-1]\cap \Z$, and let $e_t(B_M) = \E[|\xi^\Z_t \cap B_M|]$.
By translation invariance, $e_t(B_M) = p_t M$, and
\beq \nn
e_t(B_M) \leq \sum_{k\in\Z} \E[|\xi^{B_M+kM}_t \cap B_M|]
         =   \sum_{k\in\Z} \E[|\xi^{B_M}_t \cap (B_M+kM)|] = \E[|\xi^{B_M}_t|] .
\eeq
Since $M-|\xi^{B_M}_t|$ is at least as large as the number of nearest neighbor pairs in
$B_M$ that have coalesced by time $t$, we may take expectation and apply Lemma
\ref{lem:hitprob} to obtain
\beq \nn
\E[|\xi^{B_M}_t|] \leq M - (M-1)\P(|\xi^{\{0,1\}}_t|=1) \leq M- (M-1)(1-\fr{C}{\sqrt t})
< 1+ M\fr{C}{\sqrt t}.
\eeq
Therefore $p_t< 1/M + C/\sqrt{t}$. Since $M$ can be arbitrarily large for any fixed $t$,
we obtain $p_t \leq C/\sqrt{t}$.

\noindent{\bf Second Proof for continuous time.}
Let $\phi^0_t$ denote the continuous time one-dimensional voter model dual to $\xi^\Z_t$ with initial
configuration $\phi^0_0(x)=0$ for $x\in\Z\backslash\{0\}$, and $\phi^0_0(0)=1$. Let
$A_t = \{x\in\Z\ |\ \phi^0_t(x)=1\}$. The process $A_t$ is then a continuous time Markov chain on
the space of finite subsets of $\Z$. $A_t$ undergoes jumps
\beqnn
A\to A\cup\{x\}, x\notin A, \hbox{ at rate } \sum_{y\in A}\P(Y=y-x), \\
A\to A - \{x\}, x\in A, \hbox{ at rate } \sum_{y\in A^c}\P(Y=y-x). \eeqnn The rate at which
$|A_t|$ increases by 1 is $\sum_{x\in A^c}\sum_{y\in A}\P(Y=y-x)$; the rate at which $|A_t|$
decreases by 1 is $\sum_{x\in A}\sum_{y\in A^c}\P(Y=y-x)$. By translation invariance, it is not
difficult to see that these two rates are the same, and by our assumption $\P(Y=0)=0$ for continuous
time random walks, the sum of the two rates is at least 2.  Therefore, $|A_t|$ is a continuous time
simple symmetric random walk with a random rate bounded below by 2, and $|A_t|$ is absorbed at 0.
$\P(|A_t|\geq 1)$ is then bounded from above by the probability that a rate 2 simple symmetric
random walk starting from 1 does not hit 0 by time $t$, which by Lemma \ref{lem:hitprob} is bounded
above by $\fr{C}{\sqrt t}$. By the duality relation $(\ref{duality})$,
$\P(|A_t|\geq 1) = \P(\phi^0_t \not\equiv 0) = \P(0\in\xi^\Z_t)$, the lemma then follows. Note that
this proof only works for continuous time.

\blem \label{lem:negativecorrelation}
For any $A\subset Z^d$, let $\xi^A_t$ be a system of discrete time coalescing random walks on $Z^d$
starting at time 0 with one walker at every site in $A$, where all the random walks have increments
distributed as some arbitrary $\Z^d$-valued random variable $Y$. Then for any pair of disjoint sets
$B, C\subset\Z^d$, and for any time $t\geq 0$,
\beq \label{negcor1}
\P(\xi^A_t \cap B \neq \emptyset, \xi^A_t\cap C \neq \emptyset) \leq
\P(\xi^A_t \cap B \neq \emptyset) \P(\xi^A_t \cap C \neq \emptyset).
\eeq
In particular, if $d=1$, $A=\Z$, and $x,y$ are any two distinct sites in $\Z$, we have
\beq \label{negcor2}
\P(x\in\xi^{\Z}_t, y\in\xi^{\Z}_t) \leq \P(x\in\xi^{\Z}_t)\P(y\in\xi^{\Z}_t).
\eeq
\elem
\noindent{\bf Proof.} The continuous time version of this lemma is due to Arratia (see Lemma 1
in~\cite{kn:A2}). Arratia's proof uses a theorem of Harris~\cite{kn:Ha}, which breaks down for
discrete time because there are transitions between states that not comparable to each other
with respect to some partial order. However, this can be easily remedied by using an induction
argument, which we present below.

We can assume $A, B, C$ are all finite sets, since otherwise we can approximate by finite sets, and
the relevant probabilities will all converge. The main tool in the proof is again the duality between
coalescing random walks and voter models. For any pair of finite disjoint sets
$B,C\subset \Z^d$, let $\phi^{B,C}_n$ (with distribution $\nu_n^{B,C}$) be a discrete time three-type
voter model on $\Z^d$ with state space $X=\{-1,0,1\}^{\Z^d}$ and initial condition $\phi^{B,C}_0 (x)=0$
if $x\in (B\cup C)^c$; $\phi^{B,C}_0 (x)=1$ if $x\in B$; $\phi^{B,C}_0 (x)=-1$ if $x\in C$.
(See Section \ref{voter models} for the dynamics of the voter model.) Note that under the metric
$\|\eta-\zeta\| = \sum_{x\in\Z^d} 2^{-|x|_{\infty}}|\eta(x)-\zeta(x)|$, $X$ is compact.

Let $E_A^+\subset X$ (resp., $E_A^-\subset X$) be the event that some site in $A$ is assigned the
value $+1$ (resp., $-1$). Then by the duality relation $(\ref{duality})$,
$\P(\xi^A_n\cap B\neq \emptyset) = \nu_n^{B,C}(E_A^+)$,
$\P(\xi^A_n\cap C\neq \emptyset) = \nu_n^{B,C}(E_A^-)$, and
$\P(\xi^A_n\cap B\neq \emptyset, \xi^A_n\cap C\neq \emptyset) =
\nu_n^{B,C}(E_A^+\cap E_A^-)$. The correlation inequality $(\ref{negcor1})$ then becomes
\beq \label{negcor3}
\nu_n^{B,C}( E_A^+ \cap E_A^-) \leq \nu_n^{B,C} (E_A^+) \nu_n^{B,C}(E_A^-).
\eeq
We can define a partial order on the state space $X$ by setting $\eta \leq \zeta \in X$ whenever
$\eta(x) \leq \zeta(x)$ for all $x\in\Z^d$.A function $f\,:\,X\to \R$ is called increasing (resp.,
decreasing) if for any $\eta \leq \zeta$, $f(\eta) \leq f(\zeta)$ (resp., $f(\eta) \geq f(\zeta)$).
An event $E$ is called increasing (resp., decreasing) if $1_E$ is an increasing (resp., decreasing)
function. Clearly, for finite $A$, $1_{E_A^+}$ is a continuous increasing function and $1_{E_A^-}$
is a continuous decreasing function. Inequality $(\ref{negcor3})$ will follow if we show that
$\nu_n^{B,C}$ has the FKG property (see, e.g,~\cite{kn:FKG, kn:L}), i.e., for any two
continuous increasing functions $f$ and $g$,
$\int fg\, d\nu_n^{B,C} \geq \int f\, d\nu_n^{B,C}\int g\,d\nu_n^{B,C}$.

We prove this by induction. For any pair of finite disjoint sets $B, C\subset\Z^d$, $\nu_0^{B,C}$
has the FKG property because the measure is concentrated at a single configuration. Observe
that $\nu_1^{B,C}$ is a product measure and therefore also has the FKG property
(this is a simple special case of the main result of~\cite{kn:FKG}).
We proceed to the induction step, which is a fairly standard argument~\cite{kn:J}.
Assume that for all disjoint finite
sets $B$ and $C$, and for all $0\leq k\leq n-1$, $\nu_k^{B,C}$ has the FKG property. Let us denote
the collection of sites in $\Z^d$ where $\phi_{n-1}^{B,C}(x)=1$ by $B_{n-1}$, and where
$\phi_{n-1}^{B,C}(x)=-1$ by $C_{n-1}$. Then for any two continuous increasing functions $f$ and $g$,
conditioning on $\phi_{n-1}^{B,C}$, we have by the Markov property,
\beqnn
&& \int fg\, d\nu_n^{B,C} = \int \int fg\, d\nu_1^{B_{n-1},C_{n-1}}\, d\nu_{n-1}^{B,C} \\
&\geq& \int\int f\, d\nu_1^{B_{n-1},C_{n-1}}
\int g\, d\nu_1^{B_{n-1},C_{n-1}}\ d\nu_{n-1}^{B,C} \\
&\geq& \int\int f\, d\nu_1^{B_{n-1},C_{n-1}} d\nu_{n-1}^{B,C}
\int\int g\, d\nu_1^{B_{n-1},C_{n-1}}\, d\nu_{n-1}^{B,C} \\
&=& \int f\,d\nu_n^{B,C} \int g\,d\nu_n^{B,C},
\eeqnn
where we have used the FKG property for both $\nu_{n-1}^{B,C}$ and for $\nu_1^{B_{n-1}, C_{n-1}}$, and
the observation that the conditional expectations $\int\!\! f d\nu_1^{B_{n-1},C_{n-1}}\!,\!
\int\!\! gd\nu_1^{B_{n-1},C_{n-1}}$ conditioned on $\phi_{n-1}^{B,C}$ are still continuous
increasing functions. Therefore $\nu_n^{B,C}$ also has the FKG property. This concludes the induction
proof, and establishes the lemma.

Recall that $\G_\d$ and $\breve\G_\d$ denote the piecewise constant version of $\X_\d$ and $\Xc_\d$.
We can extend the definition of $d(\cdot,\cdot)$ (resp., $d_\h(\cdot,\cdot)$) to path (resp., sets
of paths) that are right continuous with left limits. Then we have
\blem \label{lem:discontapprox}
Assume $\E[|Y|^3]<+\infty$, then for any $\e>0$, $\mu_\d[ d_\h(\X_\d, \G_\d)>\e]\to 0$ and
$\breve\mu_\d[d_\h(\Xc_\d,\breve\G_\d)>\e]\to 0$ as $\d\downarrow 0$.
\elem
\noindent{\bf Proof.} Let $\e>0$ be fixed. Recall $\Lambda_{L,T} = [-L,L]\times [-T,T]$.
Choose $L$ sufficiently large, such that $\Phi(+\infty, t)-\Phi(L,t) < \e$ for all $t\in\R$, and
$\Phi(+\infty,t)-\Phi(-\infty,t) <\e$ for all $|t|>L$. Then the event $\{ d_\h(\X_\d, \G_\d)>\e\}$
occurs because either the random walk increment associated with some lattice point in
$\Lambda_{2L,L}$ exceeds $\e$ (note that the random walks are on the rescaled lattice), the probability
of which is of the order $O(\d^{-3})\P[|Y|>\etil]$ and tends to 0 as $\d\to 0$ by the assumption
$\E[|Y|^3]<+\infty$; or a random walk starting from some lattice site in $[-2L,2L]^c\times [-L,L]$
lands inside or across the spatial interval $[-L,L]$ in one step, the probability of which is bounded
by $2L\d^{-2}\sum_{k\geq\tilde L}\P[|Y|\geq k]$, which also tends to 0 as $\d\downarrow 0$ by Markov
inequality and the assumption $\E[|Y|^3]<+\infty$. This establishes the lemma for the discrete time case.
We defer the proof of the lemma for the continuous time case to Section \ref{precompact}, where
we need to carry out similar estimates.

\chapter{Proof of the Main Result}
In this chapter, we first establish the almost sure pre-compactness of $\X_\d$ and $\Xc_\d$,
so that the almost sure closures of $\X_\d$ and $\Xc_\d$ are $(\h,\f_\h)$-valued random variables.
We will then proceed to verify conditions $(B_1')$, $(T_1)$, $(I_1)$ and $(E_1)$ for $\{\X_\d\}$
and $\{\Xc_\d\}$, thus establishing the main result of this thesis, Theorem \ref{teo:walkconv}.

\section{Almost sure pre-compactness of $\X_\d$ and $\Xc_\d$} \label{precompact}
\setcounter{equation}{0}

\blem\label{lem:compact}
Assume $\E[|Y|]<+\infty$, then for any $\d\in(0,1]$, the closure of $\X_\d$ and $\Xc_\d$ in
$(\Pi,d)$, which we will also denote by $\X_\d$ and $\Xc_\d$, are almost surely compact subsets of
$(\Pi,d)$.
\elem
{\bf Proof.} We prove the lemma only for $\X_1$ and $\Xc_1$, since the proof for $\X_\d$ and $\Xc_\d$
is identical. We will show that under the mapping $(\Phi,\Psi)$, $\X_1$ and $\Xc_1$ are almost surely
equicontinuous. Note that by the properties of $(\Phi,\Psi)$, this reduces to showing the almost sure
equicontinuity of $\X_1$ and $\Xc_1$ restricted to any square $\Lambda_L =[-L,L]\times [-L,L]$.
For $\X_1$, this further reduces to showing that $\X_1$ restricted to $[-L,L]\times [k,k+1]$ is
equicontinuous for any $L>0$ and $k\in\Z$. Note that $\X_1$ restricted to $[-L,L]\times[k,k+1]$
contains either line segments connecting sites in $[-L,L]\cap\Z$ at time
$k$ to sites in $\Z$ at time $k+1$, for which there are only $2L+1$ of these;
or it contains line segments connecting sites outside $[-L,L]$ at time $k$ to some
other sites in or across $[-L,L]$ at time $k+1$. The expected number of the second type of line
segments is easily seen to be finite by our assumption $\E[|Y|]<+\infty$. Therefore almost surely,
$\X_1$ restricted to $[-L,L]\times[k,k+1]$ contains only a finite number of line segments, hence
it is equicontinuous. This proves the almost sure precompactness of $\X_1$ in $(\Pi,d)$.

The proof for $\Xc_1$ is more messy. Observe that paths in $\Xc_1$ consist of either constant
position line segments or line segments connecting consecutive jump points of a random walk, and
we are only interested in line segments that intersect $\Lambda_L$. We will show that almost surely,
we can choose $L'$ sufficiently large such that line segments in $\Xc_1$ starting from points
outside $\Lambda_{L'}$ do not intersect $\Lambda_L$. Since almost surely there are only a finite
number of jump points inside $\Lambda_{L'}$, and hence only a finite number of non-constant-position
line segments intersecting $\Lambda_L$ (note that almost surely none of the line segments is constant
in time), $\Xc_1$ restricted to $\Lambda_L$ must be equicontinuous.

Let $L$ be fixed. Let $I_{L'}$ ($L'>L$) denote the event that some line segment in $\Xc_1$
starting from some point outside $\Lambda_{L'}$ intersects $\Lambda_L$. Since $I_{L'}$ is a
descreasing family of events as $L'$ increases, it suffices to show that $\breve\mu_1(I_{L'}) \to
0$ as $L'\to +\infty$.

\begin{figure}[tp] %htbp here, top, bottom, floating page htbp, tbp, htp, tp
\begin{center}
\includegraphics{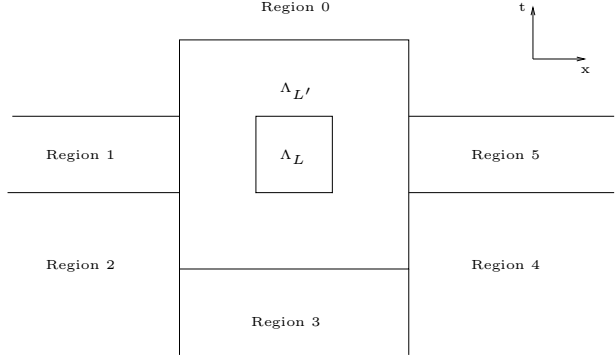}
\caption{ $\Lambda_L$ and $\Lambda_{L'}$ are centered at the origin. The solid lines divide the
complement of $\Lambda_{L'}$ into six regions. } \label{figure:Figure 2}
\end{center}
\end{figure}

We divide the complement of $\Lambda_{L'}$ into six regions as illustrated in Figure
\ref{figure:Figure 2}. Let $I_{L'}^i$ $(0\leq i\leq 5)$ denote the event that some line segment
in $\Xc_1$ starting from some point in region $i$ intersects $\Lambda_L$. Clearly
$I_{L'}^0=\emptyset$. The event $I^1_{L'}$ only occurs when the random walk increment associated
with some jump point $(x,t)$ in Region 1 exceeds $-x-L$. Since the expected number of jump
points at any site $x$ in the time interval $[-L,L]$ is $2L$, we have
\beq \label{I1}
\breve\mu_1(I^1_{L'}) \leq 2L\sum_{k\leq-L'} \P(Y\geq -L-k) = 2L\sum_{k\geq L'-L}\P(Y\geq k).
\eeq
By the assumption $\E[|Y|]<+\infty$, $\breve\mu_1(I^1_{L'})\to 0$ as $L'\to +\infty$. Similarly
$\breve\mu_1(I^5_{L'})\to 0$ as $L'\to +\infty$.

The event $I^2_{L'}$ only occurs when the random walk increment associated with some jump point
$(x,t)$ in Region 2 exceeds $-x-L$, and at the landing site there is no Poisson clock ring
during the time interval $(t, -L)$. Since the intensity of the Poisson process at each site is
1, we can estimate $\breve\mu_1(I^2_{L'})$ by
\beq \label{I2}
\breve\mu_1(I^2_{L'}) \leq \sum_{x\leq -L'}\int_{-\infty}^{-L}e^{t+L}\,\P(Y\geq-x-L)dt
= \sum_{k\geq L'-L}\P(Y\geq k).
\eeq
By the assumption $\E[|Y|]<+\infty$, $\breve\mu_1(I^2_{L'})\to 0$ as $L'\to +\infty$. Similarly
$\breve\mu_1(I^4_{L'})\to 0$ as $L'\to +\infty$. An analogous calculation shows that
$\breve\mu_1(I^3_{L'}) \leq 2L'e^{-L'+L}$, which also tends to 0 as $L'\to +\infty$. Since $I_{L'}
\subset \cup_{i=0}^5 I^i_{L'}$, the lemma then follows.

\brm
Note that the closure of $\X_\d$ $($resp., $\Xc_\d)$ in $(\Pi,d)$ is obtained from $\X_\d$
$($resp., $\Xc_\d)$ by adding all paths of the form $(f,t)$ with $t\in\d^2\Z\cup\{+\infty,-\infty\}$
$($resp., $t\in\R \cup\{+\infty,-\infty\})$ and $f\equiv +\infty$ or $f\equiv -\infty$.
\erm

We now complete the proof of Lemma \ref{lem:discontapprox}.

\noindent{\bf Proof of Lemma \ref{lem:discontapprox} for continuous time.}
In Figure \ref{figure:Figure 2}, let $(\d/\s)\Z\times\d^2\R$ be the underlying space time
lattice, and let $L$ and $L'$ be macroscopic units of space and time. Let $L$ be sufficiently
large such that $\Phi(+\infty, t)-\Phi(L,t) < \e$ for all $t\in\R$, and
$\Phi(+\infty,t)-\Phi(-\infty,t) <\e$ for all $|t|>L$. Let $L'=2L$. Keep in mind that all space
and time units are macroscopic, the event $\{ d_\h(\Xc_\d, \breve\G_\d)>\ \e \}$ only occurs when:
either for some jump point inside $\Lambda_{L'}$, the random walk increment exceeds $\e$ in magnitude
($\etil$ in unscaled units), the probability of which is of the order $\d^{-3}\P(|Y|>\etil)$ and tends
to 0 as $\d\to 0$; or for some jump point in Regions 1 and 5, the random walk penetrates the spatial
interval $[-L,L]$ in one jump; or for some jump point in Regions 2 and 4, the random walk penetrates the
spatial interval $[-L,L]$, and does not encounter another jump point until after time $-L$; or
for some jump point in Region 3, the random walk jumps and does not encounter another jump point
until after time $-L$. The probabilities of all these events can be estimated exactly as in the
proof of Lemma \ref{lem:compact} for $\Xc_1$, and they all tend to 0 as $\d\downarrow 0$. Lemma
\ref{lem:discontapprox} then follows.

%%%%%%%%%%%%%%%%%%%%%%%%%%%%%%%%%%%%%%%%%%%%%%%%%%%%%%%%%%%%%%%%%%%%%%%%%%%%%%%%%%%%%

\section{Verification of $(B_1')$} \label{B1'}
\setcounter{equation}{0}

{\bf Verification of $(B'_1).$} We first treat the discrete time case. The continuous time
case will be similar. Fix $t_0, a \in \R$, $\b>0$, $t>\b$, $\e >0$. Also fix a $\d$ and
let $\totil,\ttil, \atil$ and $\etil$ be defined from $t_0,t,a$ and $\e$ by diffusive scaling.
Then we have
\[  \mu_{\d}(|{\cal N}_{t_0,t} ([a-\e,a+\e]) | > 1)
  = \mu_1(|{\cal N}_{\totil,\ttil}
                    ([\atil-\etil,\atil+ \etil]) | > 1).
\]

If $\ttil_0 =n_0 \in\Z$, then the contribution to $\cal N$ is all due to
walkers starting from  $[\atil-\etil, \atil+\etil]\cap\Z$ at time
$n_0$. Thus we have
\begin{eqnarray}
\!\!\!&&\!\!\! \mu_1(|{\cal N}_{n_0,\ttil}
                    ([\atil-\etil,\atil+ \etil]) | > 1)
   = \mu_1(|\xi^{[\atil -\etil, \atil+\etil]\cap\Z, n_0}
         _{n_0+\ttil}| >1) \nn \\
\!\!\!&\le&\!\!\!
      \sum_{i=\lceil \atil-\etil\rceil}^{\lfloor\atil+\etil\rfloor-1}
                   \mu_1(|\xi_{n_0+\ttil}^{\{i,i+1\},n_0}|>1) \nn\\
\!\!\!&\le&\!\!\!  2\etil\mu_1(|\xi_{\ttil}^{\{0, 1\},0}|>1)
        \le 2\etil \frac{C}{\sqrt{\ttil}}
        =  \fr{2C\s\e}{\sqrt t}<\fr{2C\s\e}{\sqrt\b} \label{bineq1}
\end{eqnarray}
The first inequality follows from the observation that if the collection of
walkers starting from $[\atil -\etil, \atil+\etil]\cap\Z$ at $n_0$
has not coalesced into a single walker by $n_0+\ttil$, then there
is at least one adjacent pair of such walkers which has not coalesced by
$n_0 +\ttil$. The next inequality follows from Lemma $\ref{lem:hitprob}$.

Now suppose $\totil \in (n_0,n_0+1)$ for some $n_0\in\Z$. Note that a walker's
path can only cross $[\atil-\etil, \atil+\etil]\times\{\totil\}$ due to
the increment at time $n_0$. After the increment, at time $n_0+1$, it will either
land in $[\atil -2\etil,\atil+2\etil]$, or else outside that interval.
In the first case, the contribution of the walker's path to $\cal N$ is included in
$\xi^{[\atil -2\etil, \atil+2\etil]\cap\Z, n_0+1}_{\totil+\ttil}$, the probability
of which by our previous argument is bounded by $\fr{4C\s\e}{\sqrt t}$ times a prefactor
which approaches 1 as $\d \to 0$. In the second case, either a walker in $(-\infty, \atil+\etil]$
jumps to the right of $\atil+2\etil$, or a walker in $[\atil-\etil, +\infty)$
jumps to the left of $\atil-2\etil$, the probability of which is bounded by
\begin{eqnarray}
&&    \sum_{x=\lceil\atil-\etil\rceil}^{+\infty}\P( Y\le \atil-2\etil-x)
      + \sum_{x= -\infty}^{\lfloor\atil+\etil\rfloor}\P( Y\geq \atil+2\etil-x)
      \nn  \\
&\le &  \sum_{k=0}^{+\infty}\P(|Y|\geq k+\etil)
        \ \leq\ \sum_{k=0}^{+\infty}\fr{\E[Y^2, |Y|\geq k+\etil]}{(k+\etil)^2}
\nn\\
&\le &  \sum_{k=0}^{+\infty}\fr{\E[Y^2, |Y|\geq \etil]}{(k+\etil)^2}
        \ \leq \ \frac{2\E[Y^2, |Y|\geq\etil]}{\etil} \leq \fr{2\s^2}{\etil}
        = \frac{2\s\d}{\e}.
       \label{bineq2}
\end{eqnarray}
The next to last inequality in $(\ref{bineq2})$ is valid if we take $\d$ to be sufficiently
small. The bounds in $(\ref{bineq1})$ and $(\ref{bineq2})$ are independent of
$t_0, t>\beta$ and $a$. Taking the supremum over $t>\b, t_0$ and $a$, and letting
$\d \to 0^+$, we establish $(B'_1)$ for $\{\mu_\d\}$.

The continuous time case is similar to the discrete time case. We first scale $\Xc_\d$ back to
the $\Z\times\R$ lattice, and note that the contribution to
${\cal N}_{\totil,\ttil} ([\atil-\etil,\atil+ \etil])$ is only due to interpolated random walk
paths intersecting $[\atil-\etil,\atil+ \etil]$ at time $\totil$. Therefore
$\breve\mu_1(|{\cal N}_{\totil,\ttil} ([\atil-\etil,\atil+ \etil]) | > 1)$ can be estimated by the
union of three events: (i) event $A$, for some site $x\in [\atil-2\etil,\atil+ 2\etil]\cap\Z$,
there is no Poisson clock ring during the time interval $[\totil, \totil+\ttil]$; (ii)
event $B$, the set of coalescing random walks starting from $[\atil-2\etil,\atil+ 2\etil]\cap\Z$
at time $\totil$ has not coalesced into a single walker by time $\totil+\ttil$; (iii) event
$C$, some interpolated random walk paths intersects $[\atil-\etil,\atil+ \etil]$ at time $t_0$,
and after the intersection does not land at a site in $[\atil-2\etil,\atil+ 2\etil]\cap\Z$.
The event $\{|{\cal N}_{\totil,\ttil} ([\atil-\etil,\atil+ \etil]) | > 1\}$ is contained in
the union of the events $A$, $B$ and $C$. $\breve\mu_1(A)$ is bounded by
$4\etil e^{-\ttil}$ and tends to 0 as $\d\to 0$. $\breve\mu_1(B)$
can be estimated exactly as the computation in $(\ref{bineq1})$, and we obtain the desired
factor of $\e$ in the limit as $\d\to 0$. The event $C$ plays the same role as the event whose
probability was estimated in $(\ref{bineq2})$, and $\breve\mu_1(C)\to 0$ as $\d\to 0$,
but we will defer its proof to the next section on tightness, where we need to estimate a
similar, but more general quantity (see the paragraph above Remark \ref{rmktight}).

\bcor \label{cor:uniquepath}
Assume $\X$ (with distribution $\mu$) is a subsequential limit of $\X_\d$ (or $\Xc_\d$), then for
any deterministic point $y\in\R^2$, $\X$ has almost surely at most one path starting from $y$.
\ecor
{\bf Proof.}
It was shown in the proof of Theorem 5.3 in~\cite{kn:FINR2} that $(B'_1)$ implies
\beqnn
(B''_1)&&\!\!\!\!\!\!\!\!  \forall\b >0,
\sup_{t >\b}  \sup_{t_0,a \in \R} \mu
(|{\cal N}_{t_0,t} ([a-\e,a+\e]) | > 1) \to 0 \hbox{ as } \e
\to 0^+,
\eeqnn
the corollary then follows.

\brm
Note that if $Z^{A_\d}_\d$ is the process of coalescing random walks (either discrete or continuous
time) starting from a subset $A_\d$ of the rescaled lattice, and $Z^{A_\d}_\d$ converges in
distribution to a limit $Z$, then by the same argument as above, for any deterministic point
$y\in\R^2$, $Z$ has almost surely at most one path starting from $y$.
\erm

%%%%%%%%%%%%%%%%%%%%%%%%%%%%%%%%%%%%%%%%%%%%%%%%%%%%%%%%%%%%%%%%%%%%%%%%%%%%%%%%%%%

\section{Verification of $(T_1)$}
\setcounter{equation}{0}

In this section, we verify condition $(T_1)$ for the measures $\{\mu_\d\}$ and $\{\breve\mu_\d\}$
under the assumption $\E[|Y|^5]<+\infty$. At the end of this section, we will also show that for
$\{\X^{0_T}_\d\}$ and $\{\Xc^{0_T}_\d\}$, the set of interpolated coalescing random walk paths
starting with one walker at every site in the rescaled lattice at time 0, $\E[|Y|^3]<+\infty$ will
be sufficient to guarantee tightness.

Define $A^+_{t,u}(x_0,t_0)$ to be the event that $K$ contains a path touching both
$R(x_0,t_0;u,t)$ and (at a later time) the right boundary of the bigger rectangle $R(x_0,t_0;17u,2t)$,
and similarly define the event $A^-_{t,u}(x_0,t_0)$ corresponding to the left boundary of
the bigger rectangle. Then $A=(A^+\cup A^-)$, and writing $(T_1)$ in terms of $\mu_1$,
we argue that it is sufficient to prove
\beqnn
(T^+_1)&& \tilde g(t,u;L,T)\equiv t^{-1}\limsup_{\d\to 0^+}\,\,
            \mu_1(A^+_{\ttil,\util}(0,0))\to0\mbox{ as }t\to 0^+.
\eeqnn
The sup over $x_0, t_0$ has been safely omitted because $\mu_1$ is invariant
under translation by integer units of space and time. When $\xotil, \totil \notin\Z$,
we can bound the probability from above by using larger rectangles with vertices in
$\Z\times\Z$ and base centered at $(0,0)$. Since the argument establishing the analogous
tightness condition $(T_1^-)$ for the event $A^-$ is identical to that for $(T_1^+)$, $(T_1)$
for the measures $\{\mu_\d\}$ follows from $(T_1^+)$. Simlarly, $(T_1)$ for the measures
$\{\breve\mu_\d\}$ follows from $(T_1^+)$ with $\mu_1$ replaced by $\breve\mu_1$.

\begin{figure}[tp] %htbp here, top, bottom, floating page htbp, tbp, htp, tp
\begin{center}
\includegraphics{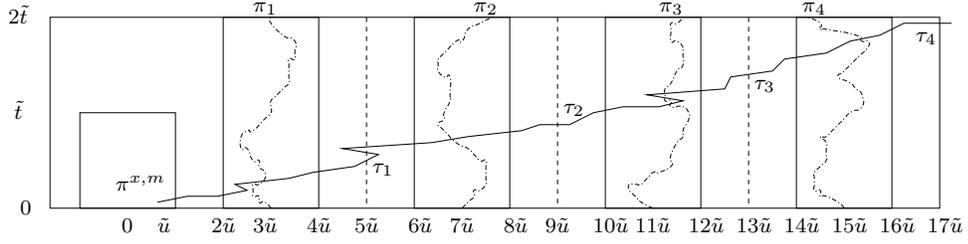}

\caption{The random walks $\pi_1,\pi_2, \pi_3$ and $\pi_4$ start from
$3\util, 7\util, 11\util$ and $15\util$ at time 0 and each stays within a distance of
$\util$ from its initial position. The random walk $\pi^{x,m}$ starts from $(x, m)$ inside the
rectangle $R(\util,\ttil)$ and exits the right boundary of the rectangle $R(17\util,2\ttil)$ at
time $\tau_4$ without coalescing with $\pi_1, \pi_2,\pi_3$ and $\pi_4$ on the way.}
\label{figure:Figure 1}
\end{center}
\end{figure}

Before we prove $(T_1^+)$, and hence $(T_1)$ for $\mu_1$ and $\breve\mu_1$, we introduce some
simplifying notation. We will abbreviate $A^+_{t,u}(0,0)$ by $A^+_{t,u}$, or just $A^+$, and
abbreviate $R(0,0;u,t)$ by $R(u,t)$. Denote the random walks (either discrete or continuous time)
starting at time 0 from $x_1 =\lceil 3\util \rceil$, $x_2=\lceil 7\util\rceil$,
$x_3=\lceil 11\util\rceil, x_4=\lceil 15\util\rceil$ by $\pi_1, \pi_2, \pi_3, \pi_4$ (with their
paths taken to be the piecewise constant version). Denote the event that $\pi_i$, $(i=1,2,3,4)$
stays within a distance $\util$ of $x_i$ up to time $2\ttil$ by $B_i$
(see Figure \ref{figure:Figure 1}). For a random walk starting from $(x,m)\in R(\util,\ttil)$,
denote the stopping times when the walker's path $\pi^{x,m}(s)$ first exceeds $5\util, 9\util, 13\util$
and $17\util$ by $\tau_1^{x,m}, \tau_2^{x,m}, \tau_3^{x,m}$ and $\tau_4^{x,m}$. We also
define $\tau_0^{x,m}=m$, and $\tau_5^{x,m}=2\ttil$. Denote the
event that $\pi^{x,m}$ does not coalesce with $\pi_i$ before time $2\ttil$ by $C_i(x,m)$. As we
shall see, the reason for choosing four paths $\pi_i$ is because each path contributes a factor of
$\d$ to our estimate of the $\mu_1$ (resp., $\breve\mu_1$) probability in $(T_1^+)$, and an overall
factor of $\d^4$ is needed to outweigh the $O(\d^{-3})$ number of lattice points (resp., jump points)
in the rectangle $R(\util,\ttil)$ from where a random walk can start. We are now ready to prove
$(T_1^+)$ for the discrete time case $\mu_1$. The proof of $(T_1^+)$ for the continuous time case
$\breve\mu_1$ is similar and will be discussed afterwards.

{\bf Verification of $(T_1^+)$ for $\mu_1$.} First we can assume $\ttil\in\Z$, since we can always
replace $\ttil$ by $\lceil\ttil\rceil$ which only enlarges the event
$A^+$. The contribution to the event $A^+$ is either due to random walk
paths that originate from within $R(\util, \ttil)$, or paths that cross
$R(\util,\ttil)$ without landing inside it after the crossing.
Denote the latter event by $D(\util,\ttil)$. Then
\begin{eqnarray}
\!\!\!\!\!\!\mu_1(A^+_{\ttil,\util})\!\!\!&\leq&\!\!\!
\mu_1(\bigcup_{i=1}^4B^c_i) + \mu_1(D(\util,\ttil)) \nn\\
\!\!\!\!\!\!&+&\!\!\!  \mu_1(\bigcap_{i=1}^4B_i;\
            \exists (x,m)\in R(\util,\ttil)
             \ s.t.\ \bigcap_{i=1}^4C_i(x,m), \tau_4^{x,m}<2\ttil).  \label{tineq1}
\end{eqnarray}
By Lemma $\ref{lem:bar}$, the first term on the right hand side of $(\ref{tineq1})$ is
of order $o(t)$ after taking the limit $\d\to 0$. For the second
event in $(\ref{tineq1})$ to occur, either a walker at a site in
$(-\infty,-\util] \times \{n\}$ jumps to a site in $[\util,+\infty)$
in one step, or a walker in $[\util,+\infty)\times\{n\}$ jumps to a site in
$(-\infty,-\util]$ in one step for some $n\in [0,\ttil-1]\cap\Z$. Denote
the event just described by $D'(\util,n)$, then
$\mu_1(D(\util,\ttil)) \leq \sum_{n=0}^{\ttil-1}\mu_1(D'(\util,n))$. From this we
see that repeating the calculations in $(\ref{bineq2})$ for $(B'_1)$ under
the assumption $\E[|Y|^3]<+\infty$ will guarantee that
$\mu_1(D(\util,\ttil))\to 0$ as $\d\to 0^+$.

To estimate the third term in $(\ref{tineq1})$ (see Figure \ref{figure:Figure 1} for
an illustration of the event), we first treat the case
of a fixed $(x,m) \in R(\util,\ttil)$. Suppressing $(x,m)$ from $\pi^{x,m}$,
$C_i(x,m)$ and $\tau^{x,m}_i$, we have
\begin{eqnarray}
\!\!\!\!\!\!\!\!\!\!\!\!&&\!\!\! \mu_1\{ {\rm for}\ (x,m)\ {\rm fixed},\ \bigcap_{i=1}^4B_i,\
             \bigcap_{i=1}^4C_i,\ \tau_4<2\ttil\}  \nn\\
\!\!\!\!\!\!\!\!\!\!\!\!&\leq&\!\!\! \mu_1\{\pi(\tau_1)>5\fr{1}{2}\util\
            \ or\ \pi(\tau_2)>9\fr{1}{2}\util\ \ or\
             \pi(\tau_3)>13\fr{1}{2}\util\} \nn\\
\!\!\!\!\!\!\!\!\!\!\!\!&+&\!\!\!\!    \mu_1\{\pi(\tau_1)\leq 5\fr{1}{2}\util,
        \pi(\tau_2)\leq 9\fr{1}{2}\util,
                \pi(\tau_3)\leq 13\fr{1}{2}\util, \tau_4<2\ttil,
        \bigcap B_i, \bigcap C_i\}. \label{tineq2}
\end{eqnarray}
The first part is bounded by
\begin{eqnarray}
 3\sup_{x\in\Z^-} \P_x[\pi^x(\tau_{0^+})>\fr{\util}{2}]
 &\leq& 3\, (\fr{2}{\util})^3
         \sup_{x\in\Z^-}\E_x[(\pi^x(\tau_{0^+}))^3, \pi^x(\tau_{0^+})>\fr{\util}{2}]
                        \nn \\
 &\leq& \fr{24}{u^3\s^3}\,\d^3 \omega(\d),  \label{tineq3}
\end{eqnarray}
where $\omega(\d)\to 0$ as $\d\to 0$. The last inequality is due to the uniform integrability
of the third moment of the overshoot distribution, which follows from our assumption
$\E[Y^5]<+\infty$ and Lemma $\ref{lem:overshoot}$.

For the second $\mu_1$ probability in $(\ref{tineq2})$, denote the event that none of the
conditions listed are violated by time $t$ by $G_t$. If $\tau_1>t$,
we interpret an inequality like $\pi(\tau_1)\leq 5\fr{1}{2}\util$ as not having been violated
by time $t$. $G_t$ is then a nested family of events, and the second probability in
$(\ref{tineq2})$ becomes
\begin{eqnarray*}
\mu_1(G_{2\ttil}) = \mu_1(G_{\tau_5}) = \mu_1(G_m) \prod_{k=1}^5
            \mu_1(G_{\tau_k} | G_{\tau_{k-1}})
          < \prod_{k=1}^4 \mu_1(G_{\tau_k} | G_{\tau_{k-1}}).
\end{eqnarray*}
Denote the history of the random walks
$\pi^{x,m}, \pi^1, \pi^2,\pi^3$ and $\pi^4$ up to time $t$ by $\Pi_t$, and denote
expectation with respect to the conditional distribution of $\Pi_t$ conditioned on
the event $G_t$ by $\E_t$. Then for $k=1,2,3,4$,
\beqnn
\mu_1(G_{\tau_k} | G_{\tau_{k-1}})
= \E_{\tau_{k-1}} \big[ \mu_1(G_{\tau_k} | \Pi_{\tau_{k-1}} \in G_{\tau_{k-1}}) \big],
\eeqnn
where the $\mu_1$ probability on the right hand side is conditioned on a given realization
of $\Pi_{\tau_{k-1}} \in G_{\tau_{k-1}}$, which is a positive probability event. For any
$\Pi_{\tau_{k-1}} \in G_{\tau_{k-1}}$, we have by the strong Markov property that
\begin{eqnarray}
\mu_1(G_{\tau_k} | \Pi_{\tau_{k-1}} \in G_{\tau_{k-1}}) &=& \mu_1[G_{\tau_k} | \pi^{x,m}(\tau_{k-1}), \pi^i(\tau_{k-1}), i=1,2,3,4] \nn\\
&\leq& C(t,u) \d,  \label{tineq4}
\end{eqnarray}
where the inequality follows from Lemma $\ref{lem:crossbar}$ for $\d$ sufficiently small.
Thus $\mu_1(G_{\tau_k} | G_{\tau_{k-1}}) \leq C(t,u)\d$, and
$\mu_1(G_{2\ttil}) \leq C^4(t,u)\d^4$. We then have
\begin{eqnarray}
&&     \mu_1(\bigcap_{i=1}^4B_i;\ \exists (x,m)\in R(\util,\ttil),
             \ s.t.\ \bigcap_{i=1}^4C_i(x,m), \tau_4^{x,m}<2\ttil) \nn\\
&\leq& \sum_{x\in [-\util,\util]\cap\Z}\,\,\sum_{m\in [0,\ttil]\cap\Z}
           \mu_1\big[{\rm for}\ (x,m)\ {\rm fixed}, \bigcap_{i=1}^4B_i,\
             \bigcap_{i=1}^4C_i,\ \tau_4^{x,m}<2\ttil\big] \nn\\
&\leq& [\fr{24}{u^3\s^3}\,\d^3 \omega(\d) + C^4(t,u)\d^4]\ 2\util\ \ttil = \omega'(\d),
\label{tineq5}
\end{eqnarray}
where $2\util\ttil = O(\d^{-3})$ and hence $\omega'(\d)\to 0$ as $\d\to 0$. Thus the last
two terms in $(\ref{tineq1})$
go to 0 as $\d\to 0$, and the first term is of order $o(t)$ after taking the limit
$\d\to 0$. Together they give $(T^+_1)$ for the measure $\mu_1$.

{\bf Verification of $(T_1^+)$ for $\breve\mu_1$.} Analogous to $(\ref{tineq1})$, the event $A^+$ is
contained in the union of three events: (i) $\cup_{i=1}^4 B_i^c$, i.e., one of the four walks
$\pi_1,\cdots,\pi_4$ fails to stay within a distance $\util$ of its starting point $x_i$ before
time $2\ttil$; (ii) the event $\breve D(\util,\ttil; \fr{3}{2}\util,2\ttil)$, some
interpolated random walk path first intersects $R(\util,\ttil)$, and then lands at a jump point
outside $R(\fr{3}{2}\util, 2\ttil)$; (iii) or the event that $\pi_1,\cdots,\pi_4$ all stay
within a distance $\util$ of its starting point $x_i$ before time $2\ttil$, and some random walk
path $\pi^{x,m}$ starting from one of the jump points $(x,m)\in R(\fr{3}{2}\util, 2\ttil)$ (by
definition, there are two random walk paths starting from any jump point, here we take the random
walk path that jumps immediately) exits from the right boundary of $R(17\util, 2\ttil)$
without first coalescing with any of the $\pi_i's$.

Using the continuous time version of Lemma \ref{lem:bar}, the probability of the first event
$\cup_{i=1}^4B_i^c$ is of order $o(t)$ after taking the limit $\d\to 0$. For the third event,
note that conditioned on the existence of a jump point at $(x,t)\in R(\fr{3}{2}\util, 2\ttil)$,
by the Markov property of Poisson process, we can apply the computations in
$(\ref{tineq2})$--$(\ref{tineq4})$ to find that the probability of a random  walk starting from
$(x,t)$, jumpping immediately, and exiting the right boundary of $R(17\util,2\ttil)$ while the
events $B_i$ all hold is of order $o(\d^3)$. Since the expected number of jump points in
$R(\fr{3}{2}\util, 2\ttil)$ is of order $\d^{-3}$, the probability of the third event is of the
order $\d^{-3} o(\d^3)$, which tends to 0 as $\d\to 0$. To conclude the proof of $(T_1^+)$ for
$\breve\mu_1$, it then only remains to show that the probability of the second event,
$\breve\mu_1[\breve D(\util,\ttil; \fr{3}{2}\util,2\ttil)]\to 0$ as $\d\to 0$.

The computations to bound $\breve\mu_1[\breve D(\util,\ttil; \fr{3}{2}\util,2\ttil)]$ is essentially
the same as that for $\breve\mu_1(I_{L'})$ in our earlier proof of the almost sure precompactness
of $\breve\X_1$ in Lemma \ref{lem:compact}. In Figure \ref{figure:Figure 2}, we replace the inner square
$\Lambda_L$ by $R(\util,\ttil)$, and the outer square $\Lambda_{L'}$ by
$R(0,-\ttil;\fr{3}{2}\util,3\ttil)$. We
can assume that for all $x\in [-\fr{3}{2}\util,\fr{3}{2}\util]\cap\Z$, there is at least one
poisson clock ring during the time interval $[\ttil, 2\ttil]$, since the probability of the
complentary event tends to 0 as $\d\to 0$. Then no constant-position line segment in $\Xc_1$
can intersect $R(\util,\ttil)$ without landing at a jump point in $R(\fr{3}{2}\util,2\ttil)$.
For non-constant-position line segments in $\Xc_1$ that originate from jump points inside
$R(0,-\ttil;\fr{3}{2}\util,3\ttil)$ and intersect $R(\util,\ttil)$ without landing at jump
points in $R(\fr{3}{2}\util,2\ttil)$, the probability is bounded by the expected number of jump
points in $R(0,-\ttil;\fr{3}{2}\util,3\ttil)$, which is of order $\d^{-3}$, times the
probability that the random walk increment $Y$ has $|Y|>\util/2$. Since $\E[|Y|^3]<+\infty$,
this product tends to 0 as $\d\to 0$. To estimate the probability of having line segments in
$\Xc_1$ that originate outside $R(0,-\ttil;\fr{3}{2}\util,3\ttil)$ and intersect
$R(\util,\ttil)$, the computation is exactly the same as that for $\breve\mu_1(I_{L'})$ in our
earlier proof of the almost sure precompactness of $\Xc_1$. Assuming $\E[|Y|^3]<+\infty$, we
find that the probability of such events also tend to 0 as $\d\to 0$.

After translation in space and time, the event $C$ defined in Section \ref{B1'} at the end of
the verification of $(B'_1)$ for the continuous time case is then easily seen to be a subset
of the event $\breve D(\etil,\ttil/2; 2\etil, \ttil)$. Therefore $\breve\mu_1(C)\to 0$ as $\d\to 0$.

\brm \label{rmktight}
The only place in this thesis where we need the assumption $\E[|Y|^5]<+\infty$ is in
$(\ref{tineq3})$. We only need $\E[|Y|^3]<+\infty$ to estimate $\mu_1(D(\util,\ttil))$
in $(\ref{tineq1})$ and to apply Lemma $\ref{lem:crossbar}$ in $(\ref{tineq4})$.
A finite third moment is the minimal moment condition for the convergence of $\X_\d$ and $\Xc_\d$.
For any $\e>0$, there are choices of $Y$ satisfying $\E[|Y|^{3-\e}]<+\infty$, but with
$\mu_1(D(\util,\ttil))\to1$ $($resp., $\breve\mu_1[\breve D(\util,\ttil; \fr{3}{2}\util,2\ttil)] \to 1)$
as $\d\to0$ for all $t>0$, which implies $\{\X_\d\}$ $($resp., $\{\Xc_\d\})$ is not tight.
\erm

Let $\X_\d^{0_T}$ with distribution $\mu_\d^{0_T}$ (resp., $\Xc_\d^{0_T}$ and $\breve\mu_\d^{0_T}$)
denote the $(\h,\f_\h)$-valued random variable consisting of interpolated discrete time
(resp., continuous time) coalescing random walk paths on the rescaled lattice starting with one walker
at every site in $(\d/\s)\Z$ at time 0. We expect tightness for $\{\X_\d^{0_T}\}$ and
$\{\Xc_\d^{0_T}\}$ to hold under much weaker moment assumptions on the random walk increment $Y$. Indeed,
\blem \label{lem:tight3m}
If $\E[|Y|^3]<+\infty$, then $\{\X_\d^{0_T}\}$ $($resp., $\{\Xc_\d^{0_T}\})$ form a tight family of
$(\h, \f_\h)$ valued random variables.
\elem

\noindent{\bf Proof.} We only prove the lemma for $\{\X_\d^{0_T}\}$, the proof for $\{\Xc_\d^{0_T}\}$
is analogous. Let $\Lambda_{L,T} = [-L,L]\times [-T,T]$. Let $F_{u,t;L,T}$ denote the
event (in $\f_\h$) that $K$ (in $\h$) contains a path $(f, t_0)$ with $(f(t_1), t_1)\in\Lambda_{L,T}$
and $|f(t_2)-f(t_1)|\geq u$ for some $t_0\leq t_1 < t_2\leq t_1+t$. Recalling the arguments leading to
the formulation of the tightness condition $(T_1)$ in~\cite{kn:FINR2}, a sufficient
condition for the family of measures $\{\mu_\d^{0_T}\}$ on $(\h, \f_\h)$ to be tight is that,
\beqnn
(T_0)\ \  \hbox{ For any } u>0, \, L,T>\!\!> u,\,
\limsup_{\d\downarrow 0} \mu_\d^{0_T}(F_{u,t;L,T}) \to 0 \hbox{ as } t\downarrow 0.
\eeqnn
If $\X_\d^{0_T}$ contains a path $(f,0)$ with $(f(t_1),t_1)\in\Lambda_{L,T}$ and
$|f(t_1)-f(t_2)|\geq u$ for some $0\leq t_1<t_2\leq t_1+t$,
then $t_1\in [mt, mt+t)$ for some nonnegative integer $m$.
By examining the locations of the path $(f,0)$ at time $mt$ and $(m+1)t$, we see that there exists
a nonnegative integer $m_0$ (either $m$ or $m+1$) and a time $m_0 t< t'\leq (m_0+1)t$ (either
$t_1$, $t_2$ or $(m+1)t$), such that either (1) $|f(m_0 t)|\leq 2L$ and $|f(t')-f(m_0t)|\geq u/4$;
or (2) $|f(m_0 t)|>2L$ and $|f(t')|\leq L$. We will call the events that $\X_\d^{0_T}\in F_{u,t;L,T}$
and $\X_\d^{0_T}$ contains a path $(f,0)$ satisfying either condition (1) or condition (2) respectively
event (1) and event (2). Then the event $\{ \X_\d^{0_T} \in F_{u,t;L,T} \}$ is a subset of the union
of events (1) and (2).

Let $L_D = \{ ku/8\ |\ k\in\Z, |k|\leq \lceil \fr{2L}{u/8}\rceil \}$ and
$T_D = \{ mt \ |\  m\in\Z, 0\leq m\leq \lceil \fr{T}{t}\rceil \}$. Let $\bar A_{t,u}(x_0,t_0)$
denote the event (in $\f_\h$) that $K$ (in $\h$) contains a path touching the bottom of
$R(x_0,t_0;u/16,t)$ and the left or right boundary of $R(x_0,t_0; u/8,2t)$. Then event (1) is
a subset of $\cup_{(x_0,t_0)\in L_D\times T_D} \bar A_{t,u}(x_0,t_0)$. By the same
argument as in the verification of $(T_1)$ for $\mu_\d$, we have
\beqnn
(\bar T_1)\ \ \
\lim_{t\downarrow 0}\ \fr{1}{t}\ \limsup_{\d\downarrow 0} \sup_{(x_0,t_0)\in {L_D\times T_D}}
\mu_\d^{0_T} (\bar A_{t,u}(x_0,t_0)) \to 0,
\eeqnn
which also holds under the assumption $\E[|Y|^3]<+\infty$. This is
because in the verification of $(T_1)$, $\E[|Y|^5]<+\infty$ is used in $(\ref{tineq3})$ to
guarantee the random walk overshoot has finite third moment, which is then used in a Markov
inequality to outweigh the $O(\d^{-3})$ number of rescaled lattice points in $R(x_0,t_0;u,t)$.
For the event $\bar A_{t,u}(x_0,t_0)$, we are only concerned with random walks starting at
the bottom of $R(x_0,t_0; u/16,t)$, which contains $O(\d^{-1})$ number of rescaled lattice points.
Therefore we only need finite first moment for the random walk overshoot, which translates into
finite third moment for $Y$.

By $(\bar T_1)$,
\beqnn
&&\limsup_{\d\downarrow 0} \mu_\d^{0_T} [\ event\ (1)\ ]  \\
&\leq& \limsup_{\d\downarrow 0}
\mu_\d^{0_T} [\cup_{(x_0,t_0)\in L_D\times T_D} \bar A_{t,u}(x_0,t_0)] \\
&\leq& (\lceil \fr{4L}{u/8}\rceil +1)(\lceil \fr{T}{t}\rceil+1) \limsup_{\d\downarrow 0}
\sup_{(x_0,t_0)\in{L_D\times T_D}} \mu_\d^{0_T} (\bar A_{t,u}(x_0,t_0))\ ,
\eeqnn
which tends to 0 as $t\downarrow 0$. On the other hand, recall the notation $\xi^B_s$ for
a system of coalescing random walks on $\Z\times\Z$ starting with one walker at every site in
$B\subset\Z$ at time 0, we have
\beqnn
&& \limsup_{\d\downarrow 0} \mu_\d^{0_T}[\ event\ (2)\ ] \\
&\leq& 2(\lceil\fr{T}{t}\rceil +1) \limsup_{\d\downarrow 0}
\P\{ \xi^{[\tilde L,+\infty)\cap\Z}_s \cap (-\infty, 0]
\neq \emptyset \hbox{ for some } s\in [0,\ttil]\,\} \\
&\leq&  \fr{2}{\a} (\lceil\fr{T}{t}\rceil +1) \limsup_{\d\downarrow 0}
\P\{ \xi^{[\tilde L,+\infty)\cap\Z}_{\ttil} \cap (-\infty, 0]
\neq \emptyset \},
\eeqnn
where $\a$ is some positive constant depending only on the random walk increment $Y$.
Observe that a nondegenerate random walk with mean zero and finite variance starting at 0 will
at any later time have a minimal probability $\a>0$ (independent of time) of being on the negative axis.
If we condition on the time and location when some walker in $\xi^{[\tilde L,+\infty)\cap\Z}_s$
first reaches $(-\infty,0]$, then the second inequality is immediate.
By the duality and the natural coupling between coalescing random walks and voter models (see
Section \ref{voter models}),
the event $\{\xi^{[\tilde L,+\infty)\cap\Z}_{\ttil} \cap (-\infty, 0] \neq \emptyset  \}$ is equivalent
to the event $\{\phi^{\Z^-}_{\ttil}(x) = 1 \hbox{ for some } x\in [\tilde L,+\infty)\cap\Z\}$ for
the dual voter model $\phi^{\Z^-}_s$ with initial condition $\phi^{\Z^-}_0(x) = 1$ if $x\in\Z^-\cup\{0\}$
and $\phi^{\Z^-}_0(x)=0$ if $x\in\Z^+$; which is also equivalent to the event that the right
boundary of the corresponding voter model interface $r_s$ satisfies $r_{\ttil}\geq \tilde L$ at time
$\ttil$ (see Section \ref{interface} for more details on the voter model interface).
A result of Cox and Durrett~\cite{kn:CD} states that if the random walk increment $Y$ has finite third
moment, then $r_s/(\s\sqrt{s})$ converges in distribution to a standard Gaussian variable as
$s\to +\infty$. Therefore
\beqnn
&&\limsup_{\d\downarrow 0} \P\{ \xi^{[\tilde L,+\infty)\cap\Z}_{\ttil} \cap (-\infty, 0]
\neq \emptyset \} \\
&=& \limsup_{\d\downarrow 0}\P_{\phi^{\Z^-}_{\ttil}}(r_{\ttil} \geq \tilde L )
= \int_{\fr{L}{\sqrt t}}^{+\infty} \fr{1}{\sqrt {2\pi}} e^{-\fr{x^2}{2}} dx \
< \ e^{-\fr{L^2}{2t}}.
\eeqnn
Therefore $\lim_{t\downarrow 0} \limsup_{\d\downarrow 0} \mu_\d^{0_T}[\ event (2)\ ] =0$. Together
with our previous estimate for the event (1), this establishes $(T_0)$, and hence the lemma.

%%%%%%%%%%%%%%%%%%%%%%%%%%%%%%%%%%%%%%%%%%%%%%%%%%%%%%%%%%%%%%%%%%%%%%%%%%%%%%%%%%%%%

\section{Verification of $(I_1)$}
\setcounter{equation}{0}

Our verification of $(I_1)$ follows a similar line of argument as in the paper of Ferrari,
Fontes and Wu~\cite{kn:FFW}. We define three sets of random walks: $\{\pi^i_\d\}_{1\leq i\leq m}$,
a family of $m$ independent random walks on the rescaled lattice $(\d/\s)\Z\times \d^2\Z$
($(\d/\s)\Z\times \d^2\R$ for continuous time); $\{\pi^i_{\d,f}\}_{1\leq i\leq m}$, the family of
$m$ coalescing random walks constructed from $\{\pi^i_\d\}$ by applying a mapping $f$ to
$\{\pi^i_\d\}$ such that two walks coalesce as soon as their paths coincide (recall that $\pi^i_\d$
denote the piecewise constant version of the random walk path);
and $\{\pi^i_{\d,g}\}_{1\leq i\leq m}$, an auxiliary family of $m$ coalescing walks
constructed by applying a mapping $g$ to $\{\pi^i_\d\}$ such that two walks coalesce as soon as their
paths cross (i.e., coincide or interchange relative order; note that random walks in $\{\pi^i_{\d,g}\}$
coalesce earlier than they do in $\{\pi^i_{\d,f}\}$). Here $\{\pi^i_\d\},\{\pi^i_{\d,f}\}$ and $\{\pi^i_{\d,g}\}$
all denote the piecewise constant version of the random walk paths. We will denote their linearly
interpolated counterpart by $\{\k^i_\d\},\{\k^i_{\d,f}\}$ and $\{\k^i_{\d,g}\}$.
If we pretend for the moment that weak convergence makes sense for piecewise constant paths without
resorting to Skorohod topology, then by Donsker's invariance principle, $\{\pi^i_\d\}$ ``converge
weakly'' to a family of independent Brownian motions $\{{\cal B}^i\}_{1\leq i\leq m}$. As we will see,
the mapping $g$ is almost surely continuous with respect to $\{{\cal B}^i\}$, and
$\{{\cal B}^i_g\}_{1\leq i\leq m}$ is distributed as coalescing Brownian motions. Therefore by the
Continuous Mapping Theorem for weak convergence, $\{\pi^i_{\d,g}\}$ ``converge weakly'' to the
coalescing Brownian motions $\{{\cal B}^i_g\}$. Finally to show that $\{\k^i_{\d,f}\}$ also converges
weakly to $\{{\cal B}^i_g\}$, we will prove that the distance between the two versions of coalescing
walks $\{\pi^i_{\d,f}\}$ and $\{\pi^i_{\d,g}\}$ converges to 0 in probability, and the distance between
the linearly interpolated version $\{\k^i_{\d,f}\}$ and the piecewise constant version
$\{\pi^i_{\d,f}\}$ also converges to 0 in probability.

We introduce more notation. Let $\cd$ be any deterministic countable dense subset of
$\R^2$. Let  $y^1=(x^1,t^1),\ldots,y^m = (x^m,t^m)\in\cd$ be fixed, and let
$\B^1,...,\B^m$ be independent Brownian motions starting from
$y^1,...,y^m$. For a fixed $\d$, denote
$\lceil \ytil^i\rceil = (\lceil\xtil^i\rceil, \lceil\ttil^i\rceil)$
(resp., $\lceil \ytil^i\rceil = (\lceil\xtil^i\rceil, \ttil^i)$ for the continuous time case),
where $\xtil^i=\s\d^{-1}x^i$ and $\ttil=\d^{-2}t^i$ as defined in Chapter \ref{randomwalk},
and let $y^i_\d$ denote $\lceil \ytil^i\rceil$'s space-time position after diffusive
scaling on the rescaled lattice $(\d/\s)\Z\times\d^2\Z$ (resp., $(\d/\s)\Z\times\d^2\R$).
Let $\pitil^i$ ($i=1,\cdots,m$) be independent random walks in the $\Z\times\Z$
(resp., $\Z\times\R$) lattice starting from $\lceil\ytil^i\rceil$. We regard $(\B^1,...,\B^m)$,
and $(\pitil^1,...,\pitil^m)$ as random variables in the product metric space
$(\Pi^m, d^{*m})$, where
\begin{eqnarray}
d^{*m}[(\xi^1, \ldots,\xi^m), (\zeta^1,\ldots,\zeta^m)]
    = \max_{1\leq i\leq m}\ d(\xi^i,\zeta^i)     \label{d*}
\end{eqnarray}
and $d$ is defined in $(\ref{d})$; thus $d^{*m}$ gives the product topology on $\Pi^m$.
We will also need the metric
\begin{eqnarray}
\dbar((f_1,t_1),(f_2,t_2)) = \sup_t|\hat f_1(t) - \hat f_2(t)|\vee|t_1-t_2| \label{dbar}
\end{eqnarray}
and $\dbar^{*m}$ is defined in a similar way as $d^{*m}$. If we denote the space of
paths that are right continuous with left limits by $\bar\Pi$, and let $\bar\Pi^m$ be the
product space, then $d$, $d^{*m}$, $\dbar$ and $\dbar^{*m}$ are still well defined
metric on $\bar\Pi$ and $\bar\Pi^m$.

We now define a mapping $g$ from $(\bar\Pi^m,d^{*m})$ to $(\bar\Pi^m,d^{*m})$ that constructs
coalescing paths from independent paths. The construction is such that when two paths first
cross (i.e., coincide or interchange relative order), the path with the higher index will be replaced
by the path with the lower index after the time of intersection or order exchange. This procedure
is then iterated until no more intersections take place. To be explicit, we give the following
algorithmic construction.

Let $(\xi^1, \ldots,\xi^m) \in \bar\Pi^m$, and let $T_g^{i,j}$ denote the time  when the two paths
$\xi^i$ and $\xi^j$ first intersect or interchange relative order. We start with equivalence relations
on the set $\{1,\ldots,m\}$ by setting $i\nsim j\ \forall\ i\neq j$. We then define the one step
iteration $\G$ on $(\xi^1,\ldots,\xi^m)$ and the equivalence relations by
\beq
\tau_g = \min_{1\leq i,j\leq m, i\nsim j}\ T_g^{i,j}        \label{5A}
\eeq
%\begin{eqnarray}
\beq
i^* =  \min \{j\ |\ j\sim i; \mbox{ or } j\nsim i, \ T_g^{i,j} = \tau_g\} \label{5B}
\eeq
%\left\{ \begin{array}{ll}
%         i     & \mbox{if $\forall\ j\nsim i, \ T_g^{i,j} >\tau_g$}, \\
%
%\} & \mbox{ otherwise }.
%         \end{array} \right.
%\end{eqnarray}
\begin{eqnarray}
\G\xi^i(t)  = \left\{ \begin{array}{ll}
         \xi^i(t)     & \mbox{if $t<\tau_g$}, \\
         \xi^{i^*}(t) & \mbox{if $t\geq\tau_g$}, \end{array} \right.  \label{5C}
\end{eqnarray}
and update equivalence relations by assigning $i\sim i^*$. Iterate the mapping $\G$,
and label the successive intersection times $\tau_g$ by $\tau_g^k$. Then the iteration
stops when
$\tau_g^k=+\infty$ for some $k\in \{1,2,\ldots, m\}$, i.e., either there is no more
crossing among the different equivalence classes of paths, or all the paths have
coalesced and formed a single equivalence class. Denote the final collection of paths
by $g(\xi^1,\ldots,\xi^m) = (\xi^1_g,\ldots,\xi^m_g)$. Then it's clear by the strong
Markov property, that $(\B^1_g, \ldots, \B^m_g)$ has the distribution of coalescing
Brownian motions. However $(\pitil^1_g, \ldots, \pitil^m_g)$ is not distributed
as coalescing random walks, because for nonsimple random walks, paths can cross
before the random walks actually coalesce (by being at the same space-time lattice site).

To construct coalescing random walk paths from independent random walks,
we define another mapping $f$ from $(\bar\Pi^m, d^{*m})$ to $(\bar\Pi^m,d^{*m})$ in a similar
way as we defined $g$, except in (\ref{5A})--(\ref{5C}) we replace the time of first crossing
$\tau_g$ by the time of first coincidence
\beq
\tau_f = \min_{1\leq i,j\leq m, i\nsim j}\ T_f^{i,j}, \label{5D}
\eeq
where $T_f^{i,j}$ is the first time when the two paths $\xi^i$ and $\xi^j$ coincide.
We will label the successive coincidence times by $\tau^k_f$. Also denote
$f(\xi^1,\cdots,\xi^m)$  by $(\xi^1_f,\cdots,\xi^m_f)$. It is then clear that
$(\pitil^1_f, \ldots, \pitil^m_f)$ is distributed as coalescing random walks
starting from $(\lceil\ytil^1\rceil, \ldots, \lceil\ytil^m\rceil)$ in the unscaled
lattice $\Z\times\Z$ (or $\Z\times\R$). We shall denote the diffusively rescaled versions of
$(\pitil^1,\ldots,\pitil^m)$, $(\pitil^1_g,\ldots,\pitil^m_g)$ and
$(\pitil^1_f,\ldots,\pitil^m_f)$ by $(\pi^1_\d,\ldots,\pi^m_\d)$,
$(\pi^1_{\d,g},\ldots,\pi^m_{\d,g})$ and $(\pi^1_{\d,f},\ldots,\pi^m_{\d,f})$. We
need the following lemma to prove $(I_1)$.

\blem \label{lem:together}
$\forall \ \e>0$, $\P\{ d^{*m}[(\pi^1_{\d,f}, \ldots, \pi^m_{\d,f}),
(\pi^1_{\d,g}, \ldots, \pi^m_{\d,g})] \geq \e \} \to 0$ as $\d\to 0^+$.
\elem
{\bf Proof.} From the definition of $d$ and $\dbar$ in $(\ref{d})$ and
$(\ref{dbar})$, it is clear that $d((f_1,t_1),(f_2,t_2)) \leq \dbar((f_1,t_1),(f_2,t_2))$
for any $(f_1,t_1), (f_2, t_2) \in\bar\Pi$. Therefore it is sufficient to prove the lemma with
$d^{*m}$ replaced by $\dbar^{*m}$. In terms of random walks in the unscaled lattice, the lemma can
be stated as
\beq \label{5E}
\forall\ \e>0, \P\{\dbar^{*m}[(\pitil^1_f, \ldots, \pitil^m_f),
(\pitil^1_g, \ldots, \pitil^m_g)] \geq \etil \} \to 0\ as\ \d\to 0^+.
\eeq

We first prove $(\ref{5E})$ for $m=2$. Note that for $m=2$, $\pitil^1_f=\pitil^1_g=\pitil^1$,
hence $\dbar^{*2}[(\pitil^1_f, \pitil^2_f), (\pitil^1_g, \pitil^2_g)] = \dbar\,(\pitil^2_f,\pitil^2_g)$.
Let $\tilde T_g^{1,2}$ denote the first time when $\pitil^1$ and $\pitil^2$ cross, and let
$\tilde T_f^{1,2}$ denote the first time when the two
walks coincide. Also let $l(0,n)$ denote the maximum distance over all time between two coalescing
random walk paths $\pi^{0,0}$ and $\pi^{n,0}$ starting at 0 and $n$ at time 0. Then by the strong
Markov property, and conditioning at time $\tilde T_g^{1,2}$,
\beq \label{5F}
\P [\ \dbar(\pitil^2_f, \pitil^2_g)\ \geq \etil\ ]
\leq \sum_{n=1}^{+\infty}\P[|\pitil^1(\tilde T_g^{1,2}) - \pitil^2(\tilde T_g^{1,2})| =n]\
\P[l(0,n)\geq \etil].
\eeq
The first probability in the summand converges to a limiting probability distribution as $\d\to 0$
by applying Lemma $\ref{lem:overshootlimit}$ to $(\pitil^1 -\pitil^2)$. The second probability
converges to 0 for every fixed $n$ by Lemma $\ref{lem:excursion}$. This proves $(\ref{5E})$
for $m=2$.

For $m > 2$, let $\tilde T_f^{i,j}$ and $\tilde T_g^{i,j}$ denote respectively the first
time when the two independent walks $\pitil^i$ and $\pitil^j$ coincide or interchange relative order.
As usual, let $T_{\d,f}^{i,j} = \d^2\tilde T_f^{i,j}$ and $T_{\d,g}^{i,j} = \d^2\tilde T_g^{i,j}$.
By Donsker's invariance principle, the interpolated paths $(\k^1_\d,\cdots,\k^m_\d)$
converge in distribution to $({\cal B}^1,\cdots,{\cal B}^m)$ as $(\Pi^m,d^{*m})$ valued random variables.
By Skorohod's representation theorem~\cite{kn:Bi,kn:Du}, we may assume this convergence is almost
sure, i.e., $d^{*m}[(\k^1_\d,\cdots,\k^m_\d), ({\cal B}^1,\cdots,{\cal B}^m)]\to 0$ almost surely.
Then by the properties of standard Brownian motions, we also have
$d^{*m}[(\pi^1_\d,\cdots,\pi^m_\d), ({\cal B}^1,\cdots,{\cal B}^m)]\to 0$ almost surely.
Note that the crossing times $T^{i,j}_{\d,g}$ as functions from $(\bar\Pi^m,d^{*m})$ to $\R$
are almost surely continuous with respect to $({\cal B}^1,\cdots,{\cal B}^m)$. Therefore
almost surely, $T_{\d,g}^{i,j}\to \tau^{i,j}$ as $\d\to 0$, where $\tau^{i,j}$ is the time of
first crossing between ${\cal B}^i$ and ${\cal B}^j$; and $\{T_g^{i,j}\}_{1\leq i<j\leq m}$ converge
jointly in distribution to $\{ \tau^{i,j}\}_{1\leq i<j\leq m}$. By the standard properties of Brownian
motion, $\{\tau^{i,j}\}_{1\leq i<j\leq m}$ are almost surely all distinct.
By an argument similar to $(\ref{5F})$, we also have
$\sup_{1\leq i<j\leq m} |T_{\d,f}^{i,j}- T_{\d,g}^{i,j}|\to 0$ in probability.
Note that in our definition of the mapping $g$ that constructs
$(\pitil^1_g,\cdots,\pitil^m_g)$ from $(\pitil^1,\cdots,\pitil^m)$, the successive times of crossing
$\{\tau_g^k\}_{1\leq k\leq m-1}$, are all times of first crossing between independent paths,
i.e., $\{\tau_g^k\}_{1\leq k\leq m-1} \subset \{\tilde T_g^{i,j}\}_{1\leq i<j\leq m}$.
The event in $(\ref{5E})$ can only occur due to: either (1) for some $\tau^k_g$ in the definition
of $g$, with $\tau^k_g = \tilde T^{i,j}_g$ for some $i$ and $j$,
$\tau^{k+1}_g \leq \tilde T^{i,j}_f$; or else, (2) whenever a coalescing takes place between two
paths $\pitil^i,\pitil^j$ in the mapping $g$, the same two paths will coalesce in the mapping $f$
before another coalescing takes place in the mapping $g$, and the event in $(\ref{5E})$ occurs
because for some $\tau^k_g$ with $\tau^k_g = \tilde T^{i,j}_g$, the distance between
the two paths $\pitil^i$ and $\pitil^j$ during the time interval $[\tilde T^{i,j}_g,\tilde T^{i,j}_f]$
exceeds $\etil$. The probability of the event (1) tends to 0 as $\d\to 0$ by our observations
that $\{T^{i,j}_{\d,g}\}_{1\leq i<j\leq m}$ converges jointly in distribution to
$\{\tau^{i,j}\}_{1\leq i<j\leq m}$, which are almost surely all distinct, and the fact that
$\sup_{1\leq i<j\leq m} |T_{\d,f}^{i,j}- T_{\d,g}^{i,j}|\to 0$ in probability. The probability of the
event (2) tends to 0 by our proof of $(\ref{5E})$ for $m=2$. This proves $(\ref{5E})$ and Lemma
\ref{lem:together}.
\vspace{4mm}

{\bf Verification of $\bf (I_1)$.} It is sufficient to show that for any sequence of $\d_n\downarrow 0$,
we can find a subsequence $\d'_n\downarrow 0$, such that $(\k^1_{\d'_n,f}, \ldots, \k^m_{\d'_n,f})$
converge in distribution to $(\B^1_g, ..., \B^m_g)$. As in the proof of Lemma \ref{lem:together},
Donsker's invariance principle implies $(\k^1_{\d_n}, \ldots, \k^m_{\d_n})$ converge in distribution
to $(\B^1, ..., \B^m)$, and by Skorohod's representation theorem, we may assum this convergence is
almost sure. Then almost surely, we also have
$d^{*m}[(\pi^1_{\d_n}, \ldots, \pi^m_{\d_n}), (\B^1, ..., \B^m)]\to 0 $.
Note that the mapping $g$ as a function from $(\bar\Pi^m,d^{*m})$ to $(\bar\Pi^m,d^{*m})$ is almost
surely continuous with respect to $(\B^1, ..., \B^m)$, therefore almost surely,
$d^{*m}[(\pi^1_{\d_n,g}, \ldots, \pi^m_{\d_n,g}), (\B^1_g, ..., \B^m_g)]\to 0$ as $\d_n\to 0$.
By Lemma \ref{lem:together}, we can choose a subsequence $\{\d'_n\}\subset \{\d_n\}$ such that,
almost surely as $\d'_n\to 0$,
$d^{*m}[(\pi^1_{\d'_n,f}, \ldots, \pi^m_{\d'_n,f}), (\pi^1_{\d'_n,g}, ..., \pi^m_{\d'_n,g})]\to 0 $,
which implies  that almost surely,
$d^{*m}[(\pi^1_{\d'_n,f}, \ldots, \pi^m_{\d'_n,f}), (\B^1_g, ..., \B^m_g)]\to 0 $.
By the properties of Brownian motion, it then follows that the linearly interpolated coalescing
random walk paths $(\k^1_{\d'_n,f}, \ldots, \k^m_{\d'_n,f})$ also converge in the metric $d^{*m}$
to $(\B^1_g, ..., \B^m_g)$ almost surely, thus completing the proof of $(I_1)$.

%%%%%%%%%%%%%%%%%%%%%%%%%%%%%%%%%%%%%%%%%%%%%%%%%%%%%%%%%%%%%%%%%%%%%%%%%%%%%%%%%%%%%

\section{Verification of $(E_1)$} \label{E1}
\setcounter{equation}{0}

As usual, we start with some notation. For an $(\h, \f_\h)$-valued random variable $X$,
define $X^{s^-}$ to be the subset of paths in $X$ which start before or at time $s$, and
for $s\leq t$ define $X^{s^-, t_T}$ to be the set of paths in $X^{s^-}$ {\it truncated}
before time $t$, i.e., replacing each path in $X^{s^-}$ by its restriction to time greater
than or equal to $t$. When $s=t$, we denote $X^{s^-, s_T}$ simply by $X^{s_T}$. Also let
$X(t)\subset \R$ denote the set of values at time $t$ of all paths in $X$. Note that
$\etahat_X(t_0,t;a,b)= |X^{t_0^-}(t_0+t) \cap (a,b)|$. We may sometimes abuse the notation and
use $X(t)$ also to denote the set of points $X(t)\times\{t\} \in \R^2$.

We recall here the definition of stochastic domination as given in~\cite{kn:FINR2}.
For two measures $\mu_1$ and $\mu_2$ on $(\h, \f_\h)$, $\mu_2$ is stochastically dominated
by $\mu_1$ ($\mu_2 << \mu_1$) if for any bounded increasing function $f$ on $(\h, \f_\h)$,
(i.e. $f(K)\leq f(K')$ if $K\subset K'$), $\int f d\mu_2 \leq \int f d\mu_1$. When $\mu_1,\mu_2$
are the distributions of two $(\h, \f_\h)$-valued random variables $X_1$ and $X_2$, we will also
denote the stochastic domination by $X_2<<X_1$. The first step of our proof is to reduce $(E_1)$ to
the following condition:
\vskip 0.3cm
\noindent$(E_1')$ If $Z_{t_0}$ is any subsequential limit of $\{X_n^{t_0^-}\}$ for any $t_0\in\R$,
then $\forall t,a,b\in\R$ with $t>0$ and $a<b$,
$\E[\etahat_{Z_{t_0}}(t_0,t; a,b)] \leq \E[\etahat_{\bar\W}(t_0,t;a,b)]= \fr{b-a}{\sqrt{\pi t}}$.
\vskip 0.3cm

\blem \label{lem:E1'}
Assuming $\{ X_n\}$ is a tight family of $(\h,\f_\h)$-valued random variables, then
$(E_1')$ implies $(E_1)$.
\elem
{\bf Proof.} Let $t_0\in\R, t>0$ be fixed, and let $X$ be the weak limit of any subsequence
$X_{n_i}$ with $n_i\to +\infty$. To prove the Lemma, it is sufficient to show that for any
$0<\e<t$, there is a further subsequence $n'_i$ such that $X_{n'_i}^{(t_0+\e)^-}$ converges
weakly to a limit $Z_{t_0+\e}$, and $X^{t_0^-} << Z_{t_0+\e}$. Because then we have
$\E[\etahat_X(t_0,t;a,b)] \leq \E[\etahat_{Z_{t_0+\e}}(t_0+\e,t-\e;a,b)]\leq\fr{b-a}{\sqrt{\pi(t-\e)}}$
by $(E_1')$, and letting $\e\to 0$ establishes $(E_1)$.

To prove the existence of $\{n'_i\}$, $Z_{t_0+\e}$, and the stochastic domination, we use a
coupling argument. Define $(\h\times\h, d^{*2}_\h)$-valued random variables
$W_{n_i} = (X_{n_i}, X_{n_i}^{(t_0+\e)^-})$, where $d^{*2}_\h$ is given by
\beqnn
 d^{*2}_\h[(K_1, K_2), (K'_1,K'_2)] = \max \{d_\h(K_1,K'_1), d_\h(K_2, K'_2)\}
\eeqnn
for $K_1,K_2,K'_1,K'_2 \in \h$. $(\h\times\h,d^{*2}_\h)$ is a complete separable metric
space. Since $\{X_{n_i}\}$ is tight, and $\{X_{n_i}^{(t_0+\e)^-}\}$ is almost surely a compact
subset of $\X_{n_i}$ for all $n_i$, $\{X_{n_i}^{(t_0+\e)^-}\}$ and $\{W_{n_i}\}$ are also
tight. Therefore we can choose a subsequence $n'_i$ such that $W_{n'_i}$ converges weakly
to a limit $W=(X',Z_{t_0+\e})$, where $X'$ is equally distributed with $X$ and $Z_{t_0+\e}$
is the weak limit of $X_{n'_i}^{(t_0+\e)^-}$.
 By Skorohod's representation theorem (see, e.g., \cite{kn:Bi, kn:Du}),
we may assume the convergence is almost sure. Then almost surely, any path $(f,t) \in X'$ with
$t\leq t_0$ is the limit of a sequence of paths $(f_{n'_i}, t_{n'_i}) \in X_{n'_i}$ with
$t_{n'_i}\to t$. Since $(f_{n'_i},t_{n'_i})$ is eventually in $X_{n'_i}^{(t_0+\e)^-}$, we also have
$(f,t)\in Z_{t_0+\e}$. Therefore $X^{'t_0^-} \subset Z_{t_0+\e}$ almost surely, and
$X^{'t_0^-} << Z_{t_0+\e}$. Since $X^{'t_0^-}$ is equally distributed with $X^{t_0^-}$,
the Lemma then follows.
\vspace{5 mm}

We now cast the condition $(E_1')$ in terms of our random variables $\{\X_\d\}$ and $\{\Xc_\d\}$.
Let $Z_{t_0}$ be any subsequential limit of $\X_\d^{t_0^-}$ (or $\Xc_\d^{t_0^-}$). Then
the validity of $(E_1')$ for $\{\X_\d\}$ and $\{\Xc_\d\}$ is a consequence of the following two lemmas,
which are also what one needs to establish to verify $(E_1)$ for general models other than coalescing
random walks $\{\X_\d\}$ and $\{\Xc_\d\}$.
\blem\label{lem:localfinite}
Let $Z_{t_0}(t_0+\e) \subset \R\times\{t_0+\e\}$ be the intersections of paths in $Z_{t_0}$
with the line $t=t_0+\e$. Then for any $\e>0$, $Z_{t_0}(t_0+\e)$ is almost surely locally finite.
\elem

\blem\label{lem:randombrownianpaths}
For any $\e>0$, $Z_{t_0}^{(t_0+\e)_T}$, the set of paths in $Z_{t_0}$ (which all start at time
$t\leq t_0$) truncated before time $t_0+\e$, is distributed as ${\cal B}^{Z_{t_0}(t_0+\e)}$,
i.e., coalescing Brownian motions starting from the random set $Z_{t_0}(t_0+\e)\subset\R^2$.
\elem

{\bf Verification of $(E_1')$.} Assume Lemmas \ref{lem:localfinite} and \ref{lem:randombrownianpaths}
for the moment. Since ${\cal B}^{Z_{t_0}(t_0+\e)} << \bar\W$, we have for $0<\e<t$
\beqnn
&&\E[\etahat_{Z_{t_0}}(t_0, t;a,b)] \\
&=& \E[\etahat_{Z_{t_0}^{(t_0+\e)_T}}(t_0+\e,t-\e;a,b)]
= \E[\etahat_{{\cal B}^{Z_{t_0}(t_0+\e)}}(t_0+\e,t-\e;a,b)] \\
&\leq& \E[\etahat_{\bar\W}(t_0+\e, t-\e;a,b)] = \fr{b-a}{\sqrt{\pi(t-\e)}}\ .
\eeqnn
Since $0<\e<t$ is arbitrary, letting $\e\to 0$ establishes $(E_1')$ for $\{\X_\d\}$ and
$\{\Xc_\d \}$.

Recall that $\G_\d$ and $\breve\G_\d$ denote the piecewise constant version of $\X_\d$ and $\Xc_\d$.
Lemma \ref{lem:localfinite} is a consequence of the following:
\blem \label{lem:densitybound}
$\forall \ t_0, t, a, b\in\R$ with $t>0$ and $a<b$, we have
\beqnn
\limsup_{\d\to 0^+} \E[\etahat_{\G_\d}(t_0,t;a,b)] \leq \fr{C(b-a)}{\sqrt t}
\eeqnn
for some $0<C<+\infty$ independent of $t_0,t,a$ and $b$. The same is true for $\breve\G_\d$.
\elem
{\bf Proof.} This follows directly from Lemma \ref{lem:densityupperbound}.

Before proving Lemma \ref{lem:localfinite}, we introduce one more notation.
Denote the space of compact subsets of $(\bar\R^2,\rho)$ by $({\cal P},\rho_{\cal P})$,
with $\rho_{\cal P}$ the induced Hausdorff metric, i.e., for $A_1, A_2\in \cal P$,
\begin{equation}
\label{rhoP}
\rho_{\cal P}(A_1,A_2)=\sup_{z_1\in A_1}\inf_{z_2\in A_2}\rho(z_1,z_2)\vee
                 \sup_{z_2\in A_2}\inf_{z_1\in A_1}\rho(z_1,z_2).
\end{equation}
Note that $({\cal P},\rho_{\cal P})$ is a complete separable metric space.

{\bf Proof of Lemma \ref{lem:localfinite}.} We prove the lemma only for the discrete time case,
the continuous time case being exactly the same. Let $Z_{t_0}$ be the weak limit of a sequence
$\{\X_{\d_n}^{t_0^-}\}$. Then $\X_{\d_n}^{t_0^-}(t_0+\e)$ converges weakly to $Z_{t_0}(t_0+\e)$
as $({\cal P},\rho_{\cal P})$ valued random variables. By Lemma \ref{lem:discontapprox},
$\rho_{\cal P}[\X_{\d_n}^{t_0^-}(t_0+\e), \Gamma_{\d_n}^{t_0^-}(t_0+\e)]\to 0$ in probability.
Therefore $\Gamma_{\d_n}^{t_0^-}(t_0+\e)$ also converges weakly to $Z_{t_0}(t_0+\e)$. For any finite
interval $(a,b)$, $\{K\in ({\cal P},\rho_{\cal P})\,:\,|K\cap (a,b)\times\R| \geq k\}$ is an open set in
$(\cal P, \rho_{\cal P})$ for any $k\in\N$. Therefore by the weak convergence of
$\Gamma_{\d_n}^{t_0^-}(t_0+\e)$ to $Z_{t_0}(t_0+\e)$ and Lemma \ref{lem:densitybound},
\beqnn
&& \E[\,|Z_{t_0}(t_0+\e) \cap (a,b)\times\R|\,] \\
&=& \sum_{k=1}^{+\infty} \P[\,|Z_{t_0}(t_0+\e) \cap (a,b)\times\R|\geq k\,] \\
&\leq& \sum_{k=1}^{+\infty}\liminf_{\d_n\downarrow 0}
\P[\,|\G_{\d_n}^{t_0^-}(t_0+\e)\cap (a,b)\times\R|\geq k\,]
\\
&\leq& \liminf_{\d_n\downarrow 0} \E[\,|\G_{\d_n}^{t_0^-}(t_0+\e)\cap (a,b)\times\R|\,] \\
&\leq&  C(b-a)/\sqrt{\e}.
\eeqnn
Lemma \ref{lem:localfinite} then follows.

It only remains to prove Lemma \ref{lem:randombrownianpaths}. We need one more lemma.
\blem \label{lem:randombrownianset}
Let $A_\d$ and $A$ be $({\cal P},\rho_{\cal P})$-valued random variables, where $A$ is
almost surely a locally finite set, $A_\d$ is almost surely a subset of
$(\d\Z/\s)\times(\d^2\Z)$ $($resp., $(\d\Z/\s)\times(\d^2\R))$, and $A_\d$ converges
in distribution to $A$ as $\d\to 0$. Conditioned on $A_\d$, let  $\X^{A_\d}_\d$
$($resp., $\Xc^{A_\d}_\d)$ be the process
of discrete time $($resp., continuous time$)$ coalescing random walks on $(\d\Z/\s)\times(\d^2\Z)$
$($resp., $(\d\Z/\s)\times(\d^2\R))$ starting from the point set $A_\d$. Then as $\d\to 0$,
$\X_\d^{A_\d}$ $($resp., $\Xc_\d^{A_\d})$ converges in distribution to $\B^A$, the process of
coalescing Brownian motions starting from a random point set distributed as $A$.
\elem
{\bf Proof.} We only prove the lemma for the discrete time case, the continuous time case being
exactly the same. We first treat the case where $A$ and $A_\d$ are deterministic and
$\rho_{\cal P}(A_\d,A) \to 0$ as $\d \to 0$. Note that $\{\X^{A_\d}_\d\}$ is tight since
$\X^{A_\d}_\d$ is almost surely a subset of $\X_\d$ and $\{\X_\d\}$ is tight. If $Z$ is a
subsequential limit of $\X^{A_\d}_\d$, then by $(I_1)$ and the remark following Corollary
$\ref{cor:uniquepath}$, there is $\mu_Z$ almost surely exactly one path starting from every
$y\in A$, and the finite dimensional distributions of $Z$ are those of coalescing Brownian
motions. Therefore $Z$ is equidistributed with $\B^A$, which proves the deterministic case.

For the nondeterministic case, it suffices to show $\E[f(\X^{A_\d}_\d)] \to \E[f(\B^A)]$
as $\d \to 0$ for any bounded continuous function $f$ on $(\h,d_\h)$. If we denote
$f_\d(A_\d) = \E[f(\X^{A_\d}_\d)|A_\d]$, and $f_\B(A) = \E[f(\B^A)|A]$, then
$\E[f(\X^{A_\d}_\d)] = \E[f_\d(A_\d)]$ and $\E[f(\B^A)] = \E[f_\B(A)]$. Since $A_\d$ converges
in distribution to $A$, by Skorohod's representation theorem~\cite{kn:Bi, kn:Du}, we can construct
random variables $A'_\d$ and $A'$ which are equidistributed with $A_\d$ and $A$, such
that $A'_\d(\omega)\to A'(\omega)$ in $\rho_{\cal P}$ almost surely. Then for almost every
$\omega$ in the
probability space where $A_\d'$ and $A'$ are defined, by the part of the proof already done
(for deterministic $A_\d$ and $A$), $\X_\d^{A'_\d(\omega)}$ converges in distribution to
${\cal B}^{A'(\omega)}$. Thus $f_\d(A'_\d(\omega))=\E[f(\X_\d^{A'_\d(\omega)})]$ $\to$
$f_{\cal B}(A'(\omega)) = \E[f({\cal B}^{A'(\omega)})]$ for almost every $\omega$.
By the bounded convergence theorem, $\E[f_\d(A'_\d)] \to \E[f_\B(A')]$ as
$\d \to 0$. Since $A'_\d$ and $A'$ are equidistributed with $A_\d$ and $A$, the lemma follows.

{\bf Proof of Lemma \ref{lem:randombrownianpaths}.} Let $Z_{t_0}$ be the weak limit of
$\X_{\d_n}^{t_0^-}$ (resp., $\Xc_{\d_n}^{t_0^-}$) for a sequence of $\d_n\downarrow 0$.
We first treat the discrete time case. By Skorohod's representation theorem,
we can assume the convergence is almost sure. Then almost surely,
$\rho_{\cal P}(\X_{\d_n}^{t_0^-}(t_0+\e), Z_{t_0}(t_0+\e))\to 0$,
and $d_\h(\X_{\d_n}^{t_0^-, (t_0+\e)_T}, Z_{t_0}^{(t_0+\e)_T})\to 0$.
Let $m_\d= \d^2 \lceil \totil+\etil \rceil$,
the first time on the rescaled lattice greater than or equal to $t_0+\e$. Using the fact that the
image of $Z_{t_0}^{(t_0+\e)_T}$ under $(\Phi,\Psi)$ is almost surely equicontinuous, it is not
difficult to see that $\rho_{\cal P}(\X_{\d_n}^{t_0^-}(m_{\d_n}), Z_{t_0}(t_0+\e)) \to 0$ and
$d_\h(\X_{\d_n}^{t_0^-, (m_{\d_n})_T}, Z_{t_0}^{(t_0+\e)_T}) \to 0$ almost surely. On the other
hand, $Z_{t_0}(t_0+\e)$ is almost surely locally finite by Lemma \ref{lem:localfinite}, and
$\X_{\d_n}^{t_0^-, (m_{\d_n})_T}$ is distributed as coalescing random walks on the rescaled lattice
starting from $\X_{\d_n}^{t_0^-}(m_{\d_n}) \subset (\d\Z/\s)\times(\d^2\Z)$. Therefore by Lemma
\ref{lem:randombrownianset}, $\X_{\d_n}^{t_0^-, (m_{\d_n})_T}$ converges weakly to
${\cal B}^{Z_{t_0}(t_0+\e)}$, and $Z_{t_0}^{(t_0+\e)_T}$ is equally distributed with
${\cal B}^{Z_{t_0}(t_0+\e)}$.

We now treat the continuous time case. By Lemma \ref{lem:discontapprox}, as $\d_n\downarrow 0$,
$d_\h(\Xc_{\d_n}^{t_0^-,(t_0+\e)_T}, \breve\G_{\d_n}^{t_0^-,(t_0+\e)_T })$ $\to 0$ in probability.
Let $U_{\d_n}^{t_0^-,(t_0+\e)_T}$ denote the interpolated version of
$\breve\Gamma_{\d_n}^{t_0^-,(t_0+\e)_T}$. Then by Lemma \ref{lem:discontapprox},
$d_\h(U_{\d_n}^{t_0^-,(t_0+\e)_T }, \breve\Gamma_{\d_n}^{t_0^-,(t_0+\e)_T })\to 0$ also in probability.
Therefore $d_\h(\Xc_{\d_n}^{t_0^-,(t_0+\e)_T }, U_{\d_n}^{t_0^-,(t_0+\e)_T })\to 0$ in probability,
and $U_{\d_n}^{t_0^-,(t_0+\e)_T}\to Z_{t_0}^{(t_0+\e)_T}$ in distribution.
On the other hand, $U_{\d_n}^{t_0^-,(t_0+\e)_T}$ is the paths of continuous time coalescing random
walks starting from a random initial configuration at time $t_0+\e$, and
$U_{\d_n}^{t_0^-,(t_0+\e)_T}(t_0+\e)$ converges weakly to $Z_{t_0}(t_0+\e)$, which is almost surely
locally finite by Lemma \ref{lem:localfinite}. Therefore, by Lemma \ref{lem:randombrownianset},
$U_{\d_n}^{t_0^-,(t_0+\e)_T}$ converges in distribution to ${\cal B}^{Z_{t_0}(t_0+\e)}$. Thus
$Z_{t_0}(t_0+\e)$ is equally distributed with ${\cal B}^{Z_{t_0}(t_0+\e)}$, and the lemma is
established.

\brm
The key to the proof of Lemma \ref{lem:randombrownianpaths} is to approximate
$\X_{\d_n}^{t_0^-, (t_0+\e)_T}$\!\!\!\!
$($resp., $\Xc_{\d_n}^{t_0^-, (t_0+\e)_T})$ by the paths of a
Markov process starting from a random space-time point set such that Lemma \ref{lem:randombrownianset}
can be applied. For discrete time coalescing random walks, the natural choice is
$\X_{\d_n}^{t_0^-, (m_{\d_n})_T}$, while for continuous time, the natural choice is
$U_{\d_n}^{t_0^-, (t_0+\e)_T}$.
\erm

\chapter{Further Results}
\vspace{-1cm}
\section{Convergence of $\X_\d^{0_T}(1)$ and $\Xc_\d^{0_T}(1)$}
\setcounter{equation}{0}

Let $\X_\d^{0_T}$ (resp.,$\Xc_\d^{0_T}$) denote the $(\h,\f_\h)$-valued random variable consisting
of the set of interpolated coalescing random walk paths on the rescaled lattice starting
with one walker at every site in $\d\Z/\s$ at time 0. In~\cite{kn:A0}, Arratia proved that, for
coalescing simple random walks, $\X^{0_T}_\d (1)$ as a point process on $\R$ converges in distribution
to $\bar\W^0(1)$, the point process on $\R$ generated at time 1 by coalescing Brownian motions starting
from every point on $\R$ at time 0, which is a stationary simple point process with intensity
$1/\sqrt{\pi}$. In~\cite{kn:A1}, Arratia stated the analogous result for nonsimple walks with zero
mean and finite second moment for its increment, but a proof was not given. In this section,
we give a proof for random walks whose increments have mean zero and finite third moment.

We first recall the space and topology on which point processes are defined. Let
$(\hat{\cal N}, {\cal B}_{\hat{\cal N}})$ be the space of locally finite counting
measures on $\R$, where ${\cal B}_{\hat {\cal N}}$ is the Borel $\s$-algebra generated
by the vague topology on $\hat{\cal N}$, i.e., for $\mu_n$, $\mu\in \hat{\cal N}$, $\mu_n$ converges
vaguely to $\mu$ if for any bounded continuous function $f : \R\to\R$ with bounded support,
$\int f d\mu_n \to \int f d\mu$. (For more background on random measures and the vague
topology, see~\cite{kn:DJ,kn:K}). The vague topology on $\hat{\cal N}$ can be metrized
so that $(\hat{\cal N}, {\cal B}_{\hat{\cal N}})$ is a complete separable metric space.
Recall that $\breve\G_\d^{0_T}$ denotes the piecewise constant version of $\Xc_\d^{0_T}$,
our result on the convergence of point processes is then the following.

\bteo \label{teo:pointprocess}
If $\E[|Y|^3]<+\infty$, then $\X_\d^{0_T}(1)$, $\G_\d^{0_T}(1)$, $\Xc_\d^{0_T}(1)$ and
$\breve\G_\d^{0_T}(1)$  as $(\hat{\cal N}, {\cal B}_{\hat{\cal N}})$-valued random variables
converge weakly to $\bar\W^0(1)$ as $\d\to 0$.
\eteo

Before we prove Theorem \ref{teo:pointprocess}, we need the following lemma.
\blem \label{lem:brownianline}
If the random walk increment satisfies $\E[|Y|^3]<+\infty$, then $\X_\d^{0_T}$ and $\Xc_\d^{0_T}$
converge in distribution to $\bar\W^0$, the subset of paths in $\bar\W$ starting at time 0.
\elem
{\bf Proof.} Note that for any countable dense subset $\cd^0 \subset \R\times\{0\}$, $\bar\W(\cd^0)$,
the closure in $(\Pi,d)$ of coalescing Brownian motion paths starting from $\cd^0$ is
equidistributed with $\bar\W^0$ by properties of the Brownian Web~\cite{kn:FINR2}.
As in the case of the convergence of $\X_\d$ and $\Xc_\d$ to $\bar\W$, we need to establish tightness,
and verify conditions $(I_1)$, $(B_1')$ and $(E_1)$, where in $(I_1)$, the countable dense
set $\cd\subset \R^2$ is now replaced by $\cd^0$, and in $(B'_1)$ and $(E_1)$, $t_0$ is set to
0. Tightness follows from Lemma \ref{lem:tight3m}. The other conditions follow
directly from their verification for $\X_\d$ and $\Xc_\d$ , which all require at most finite
third moment of $Y$.

{\bf Proof of Theorem \ref{teo:pointprocess}.} We only prove the theorem for $\X_\d^{0_T}(1)$,
the proof for $\G_\d^{0_T}(1)$, $\Xc_\d^{0_T}(1)$ and $\breve\G_\d^{0_T}(1)$ is analogous. Since
$\bar\W^0(1)$ is a simple point process, to prove
weak convergence of $\X_\d^{0_T}(1)$ to $\bar\W^0(1)$, it is sufficient to show: (i) tightness;
(ii) any subsequential limit of $\X_\d^{0_T}(1)$ is a simple point process; (iii) convergence of
the avoidance (zero) functions, i.e., for any disjoint union of a finite number of finite
intervals $A = \bigcup_{i=1}^n [a_i, b_i]$,
$\P[\X_\d^{0_T}(1)\cap A =\emptyset] \to \P[\bar\W^0(1)\cap A =\emptyset]$ as $\d\to 0$.
(See, e.g., Sections 7.3 and 9.1 in~\cite{kn:DJ}).

To prove that $\{\X_\d^{0_T}(1)\}$ is tight, it is sufficient to show that for any finite
closed interval
$[a,b]$, $\sup_{0<\d<1} \E_{\X_\d^{0_T}(1)}\big[\zeta[a,b]\big] < C$, where $\zeta[a,b]$ is the
measure of $[a,b]$ with respect to an element $\zeta\in\hat{\cal N}$, and $C<+\infty$ is a constant
depending on $[a,b]$. Since $\E_{\X_\d^{0_T}(1)}[\zeta[a,b]] < \E[\etahat_{\X_\d}(0,1;a-1,b+1)]$,
tightness follows from Lemma \ref{lem:densitybound}.

Let $Z$ be a weak limit of $\{\X_{\d_n}^{0_T}(1)\}$ in $(\hat{\cal N},{\cal B}_{\hat{\cal N}})$
with distribution $\mu_Z$.
For any $\zeta\in\hat{\cal N}$, let $\zeta_{[m,n]}$ denote $\zeta$ restricted to $[m,n]$. Also
let $A_{N,i}^m = [m+(i-1)2^{-N}, m+i2^{-N}]$. Then for all $N\in\N$, $m<n\in\Z$,
\begin{eqnarray}
\!\!\!\!\!\!\!\!\!\!\!\!\mu_{Z}\{\zeta\ |\ \zeta_{[m,n]} \hbox{ not simple }\}
\!\!\!&\leq&\!\!\! \mu_{Z}(\bigcup_{i=1}^{(m-n)2^N}\{\zeta(A_{N,i}^m)\geq2 \}) \nn\\
\!\!\!&\leq&\!\!\! (m-n)2^N\sup_{a\in[m,n]}\mu_{Z}(\zeta[a,a+2^{-N}]\geq2). \label{corineq1}
\end{eqnarray}
To prove $Z$ is a simple point process, it is then sufficient to show that
\begin{eqnarray}
\label{corineq2}
\limsup_{\e\to 0^+}\ \ \fr{1}{\e} \sup_{a\in\R}\mu_{Z}\big[\zeta[a,a+\e]\geq2 \big] =0,
\end{eqnarray}
since this implies that $\zeta_{[m,n]}$ is $\mu_{Z}$ almost surely a simple counting measure by
taking $N\to+\infty$ in $(\ref{corineq1})$. Letting $m\to+\infty$ and $n\to -\infty$
then implies that $Z$ is almost surely a simple counting measure.

Note that $\mu_{Z}[\zeta[a,a+\e]\geq2]\ \leq\ \mu_{Z}[\zeta(a-\e,a+2\e)\geq2]$, and
$\{\zeta|\zeta(a-\e, a+2\e)\geq2\}$ is an open set in $(\hat{\cal N}, {\cal B}_{\hat{\cal N}})$.
By the weak convergence of $\X_{\d_n}^{0_T}(1)$ to $Z$, we have
\begin{eqnarray}
&&\mu_{Z}[\zeta(a-\e,a+2\e)\geq2] \nn\\
&\leq& \liminf_{\d_n\downarrow 0}\ \mu_{\X^{0_T}_{\d_n}(1)}[\zeta(a-\e, a+2\e)\geq2] \nn\\
&=&
\liminf_{\d_n\downarrow 0}\ \P(|\xi^\Z_{\d_n^{-2}}\cap(\atil-\etil,\atil+2\etil)|\geq2)\nn\\
&\leq&   \liminf_{\d_n\downarrow 0} \sum_{i,j = \lceil\atil-\etil\rceil, i\neq j}
                                     ^{\lfloor\atil+2\etil\rfloor}
       \P(i, j \in \xi^\Z_{\d_n^{-2}}) \nn
\end{eqnarray}
\begin{eqnarray}
&\leq&
\label{corineq3}
\liminf_{\d_n\downarrow 0}\sum_{i,j = \lceil\atil-\etil\rceil, i\neq j}
                                     ^{\lfloor\atil+2\etil\rfloor}
        \P(i\in \xi^\Z_{\d_n^{-2}})\ \P(j\in \xi^\Z_{\d_n^{-2}}) \\
&\leq&
\liminf_{\d_n\downarrow 0}\ (3\etil+1)^2\ \P [0\in \xi^\Z_{\d_n^{-2}}]^2
\ \leq\ 9\, C^2\s^2\e^2, \label{corineq4}
\end{eqnarray}
where in $(\ref{corineq3})$, we applied Lemma \ref{lem:negativecorrelation}, and in $(\ref{corineq4})$,
we applied Lemma \ref{lem:densityupperbound}. To apply Lemma \ref{lem:negativecorrelation}, we
have implicitly assumed $\d_n^{-2} \in\N$. If $\d_n^{-2} \notin\N$, then we need to approximate and
use an argument similar to the one leading to the computation in $(\ref{bineq2})$.
This establishes $(\ref{corineq2})$, thus proving any subsequential limit of $\X_{\d_n}^{0_T}(1)$ must
be a simple point process.

We now show the convergence of the avoidance functions. Let $A = \bigcup_{i=1}^n [a_i, b_i]$
be the disjoint union of a finite number of finite intervals. By Lemma \ref{lem:brownianline},
$\X_\d^{0_T}$ converges weakly to $\bar\W^0$ as $(\h, \f_\h)$-valued random variables, so by
Skorohod's representation theorem, we may assume this convergence is almost sure.
In particular, $\X^{0_T}_\d(1)$ converges almost surely to $\bar\W^0(1)$ in $\rho_{\cal P}$ as
defined in $(\ref{rhoP})$. Since $\bar\W^0(1)$ is a stationary simple point process with intensity
$1/\sqrt{\pi}$, $\P[\partial A\cap \bar\W^0(1) \neq \emptyset]=0$. It is then easy to see that
${\bf 1}_{\X^{0_T}_\d(1)\cap A=\emptyset} \to {\bf 1}_{\bar\W^0(1)\cap A=\emptyset}$ almost
surely. By the bounded convergence theorem,
$\lim_{\d\downarrow0}\P(\X^{0_T}_\d(1)\cap A=\emptyset) = \P(\bar\W^0(1)\cap A=\emptyset)$,
thus proving the the convergence of avoidance functions and the theorem.

As a corollary of Theorem \ref{teo:pointprocess}, we have
\bcor \label{cor:density}
If the random walk increment $Y$ satisfies $\E[|Y|^3]<+\infty$,
then $p_t\equiv \P(0\in\xi^\Z_t)\sim 1/(\s\sqrt{\pi t})$ as $t\to +\infty$.
\ecor
{\bf Proof.} We prove the corollary for the discrete time case, the continuous time case being
exactly the same. Let $\d_n = 1/\sqrt{n}$, and let $f_\e(x)$, for $\e>0$, be a
continuous function with support on $[-\e, 1+\e]$, with $0\leq f_\e\leq 1$, and $f_\e\equiv 1$
on $[0,1]$. By Theorem $\ref{teo:pointprocess}$,
$\lim_{n\to+\infty}\E_{\X_{\d_n}^{0_T}(1)}[\int f_\e d\zeta] = \E_{\bar\W^0(1)}[\int f_\e d\zeta]$.
Since
\beqnn
\fr{1}{\sqrt{\pi}} < \E_{\bar\W^0(1)}[\int f_\e d\zeta] < \fr{(1+2\e)}{\sqrt{\pi}},
\eeqnn
and
\beqnn
\s\sqrt{n} p_n \leq \E_{\X_{\d_n}^{0_T}(1)}[\int f_\e d\zeta] \leq (1+2\e)\s\sqrt{n} p_n.
\eeqnn
It follows that
\beqnn
\fr{1}{(1+2\e)\sqrt{\pi}} \leq
   \liminf_{n\to +\infty} \s\sqrt{n} p_n \leq \limsup_{n\to +\infty}\s\sqrt{n} p_n
\leq \fr{(1+2\e)}{\sqrt{\pi}}.
\eeqnn
Since $\e>0$ is arbitrary, letting $\e\to 0$ establishes the corollary.
\vspace{4 mm}

Let $\phi^0_t$ be a one-dimensional voter model (either discrete or continuous time) with state space
$\{0,1\}^\Z$ and initial configuration $\phi^0_0(x)=0$ for $x\in\Z \backslash\{0\}$, and
$\phi^0_0(0)=1$. The dynamics of the model is defined in Section \ref{voter models}.
Then by the duality relation (\ref{duality}), $\P(\phi^0_t\not\equiv 0)=\P(0\in\xi^\Z_t)$.
Corollary $\ref{cor:density}$ is then equivalent to
\bcor\label{cor:extinction}
Let $\phi^0_t$ be the voter model defined above. If $\E[|Y|^3]<+\infty$, then
$\P(\phi^0_t\not\equiv0)\sim 1/(\s\sqrt{\pi t})$ as $t\to +\infty$.
\ecor
\brm
Corollaries $\ref{cor:density}$ and $\ref{cor:extinction}$ partially extend a result of Bramson
and Griffeath~\cite{kn:BG}. They proved that, for continuous time coalescing {\bf simple}
random walks in $\Z^d$, $\xi^{\Z^d}_t$, and the dual voter model $\phi^{0,d}_t$ with initial
configuration all $0'$s except for a 1 at the origin,
$p_t = \P(0\in\xi^{\Z^d}_t)=\P(\phi^{0,d}_t\not\equiv 0)$ decays asymptotically as
$1/(\sqrt{\pi t})$ for $d=1$, $\log t/(\pi t)$ for $d=2$, and $1/(\g_d t)$ for
$d\geq 3$. For $d\geq 2$, their proof also works for discrete time random walks, and as
pointed out in the remark before Lemma 2 in~\cite{kn:BCL}, can be easily extended to
much more general random walks (see~\cite{kn:BCL} for more details).
\erm

\brm \label{cor:browniancorrelation}
The following correlation inequality is valid for the point process $\bar\W^0(1)$. Let $A$, $B$ be
two disjoint open sets in $\R$, and let $O_A = \{ \zeta\in{\hat{\cal N}}\,|\, \zeta(A) \geq 1\}$ and
$O_B = \{ \zeta\in{\hat{\cal N}}\,|\, \zeta(B) \geq 1\}$. Then
\beq\nn
\mu_{\bar\W^0(1)}(O_A \cap O_B) \leq \mu_{\bar\W^0(1)}(O_A)\ \mu_{\bar\W^0(1)}(O_B).
\eeq
This negative correlation inequality for $\bar\W^0(1)$ is implicit in the work of Arratia~\cite{kn:A2};
it is also a direct consequence of Lemma \ref{lem:negativecorrelation} and Theorem
\ref{teo:pointprocess}. By similar arguments, the same correlation inequality holds for point
processes generated at time 1 by coalescing Brownian motion paths starting from
any closed space-time
region strictly below time 1.
\erm

%%%%%%%%%%%%%%%%%%%%%%%%%%%%%%%%%%%%%%%%%%%%%%%%%%%%%%%%%%%%%%%%%%%%%%%%%%%%%%%%%%%%%%%%%%
\section{Voter model interface} \label{interface}
\setcounter{equation}{0}

\noindent{\bf Voter Model Interface:} Let $\phi^{\Z^-}_t$ be a one-dimensional voter model
(either discrete or continuous time) with state space $\{0,1\}^\Z$ and initial configuration
$\phi^{\Z^-}_0(x) = 1$ if $x\in \Z^-\cup\{0\}$ and $\phi^{\Z^-}_0(x) = 0$ if $x\in\Z^+$.
At any time $t\geq 0$, $\phi^{\Z^-}_t$ will contain a leftmost 0 and a rightmost 1, whose
positions we denote by $l_t$ and $r_t$. Then $\phi^{\Z^-}_t(x) = 1$ for $x< l_t$,
$\phi^{\Z^-}_t(x) = 0$ for $x>r_t$, and the configuration of $\phi^{\Z^-}_t$ between $l_t$
and $r_t$ defines what is called the interface process, $\a_t = \phi^{\Z^-}_t(x+l_t)$ with $x\in\N$,
which is a random variable taking values in
$\big\{\xi: \N\to \{0,1\}, \sum_{x\in\N}\xi(x) < +\infty\big\}$.
Cox and Durrett proved in~\cite{kn:CD} that, for the continuous time model $\phi^{\Z^-}_t$ with
the associated random walk increment $Y$ satisfying $\E[|Y|^3]<+\infty$, the interface process $\a_t$
is an irreducible positive recurrent Markov chain. Hence the size of the interface $r_t-l_t$
is of $O(1)$ as $t\to +\infty$. They also proved that $l_t/(\s\sqrt{t})$ and $r_t/(\s\sqrt{t})$
converge in distribution to standard Gaussian variables as $t\to +\infty$. Their result should also
be valid for the discrete time model $\phi^{Z^-}_n$.

The weak convergence of $\X_\d$ and $\Xc_\d$ to the Brownian Web $\bar\W$ recovers Cox and Durrett's
result under the stronger assumption that $\E[|Y|^5]<+\infty$, but it also establishes that the
time evolutions of $l_t$ and $r_t$ converge weakly to the same Brownian motion. In the following
discussion, we will let $\bar l_t, \bar r_t$, for $t\geq 0$, denote the continuous paths constructed
from $l_t, r_t, t\geq 0$ by linear interpolation the same way we construct the interpolated path of
a random walk.

\bteo\label{teo:interface}
Let $\phi^{\Z^-}_t$, $\bar l_t$ and $\bar r_t$ be as defined above. Let $l_{t,\d}$ be $l_t$ diffusively
rescaled, i.e., $l_{t,\d}=\d\s^{-1}l_{t\d^{-2}}$, and define $\bar l_{t,\d}$, $r_{t,\d}$ and
$\bar r_{t,\d}$ similary. If $\E[|Y|^5]<+\infty$, then
$\{\bar l_{t,\d}, \bar r_{t,\d}\}$ as $(\h, d_\h)$ valued random variables converge in
distribution to ${\cal B}_t^{0,0}$, a standard Brownian motion starting at the origin at time 0.
\eteo
{\bf Proof.} We first prove the theorem for the discrete time case.
Let $\hat\X_\d = \{(f(-t), t\leq-t_0)\ |\ (f(t), t\geq t_0)\in\X_\d \}$, i.e. the
coalescing random walks running backward in time on the rescaled lattice.
By the duality relation (\ref{duality}) and the natural coupling between voter models and coalescing
random walks, we have for all $x\in\Z$ and $t\in\N$, $\phi^{\Z^-}_t(x) = \phi^{\Z^-}_0(\hat\k^{x,t}_0)$
almost surely, where $\hat\k^{x,t}_0$ is the location at time 0 of the backward random walk path in
$\hat\X_1$ starting at $x$ at time $t$.

By Theorem \ref{teo:walkconv}, $\hat\X_\d$ converges weakly to $\hat{\bar\W}$ (the backward Brownian
Web) as $\d\to 0$; and by Skorohod's representation theorem, we may assume this convergence is almost
sure. $\hat{\bar\W}$ uniquely determines a dual (forward) Brownian web $\bar\W$, which is equally
distributed with the standard Brownian web. The pair $(\bar\W, \hat{\bar\W})$ forms what is called
the {\em double Brownian web} $\bar\W^D$ with the property that, almost surely, paths in
$\bar\W$ and $\hat{\bar\W}$ {\em reflect} off each other and never
cross~\cite{kn:FINR2, kn:TW, kn:STW, kn:A1}.

Let $\phi^{\Z^-}_{t,\d}$ be the diffusively rescaled voter model dual to $\hat\X_\d$ with the natural
coupling between the two models. We will overload
the notation and let $l_{t,\d}$ and $r_{t,\d}$ also denote the piecewise constant interface boundary
lines of $\phi^{\Z^-}_{t,\d}$, and let $\bar l_{t,\d}$ and $\bar r_{t,\d}$ denote their linearly
interpolated counterparts. To prove the theorem, it is then sufficient to show that almost surely,
$\bar l_{t,\d}$ and $\bar r_{t,\d}$ converge in $(\Pi, d)$ to $\B^{0,0}_t\in\bar\W$, the unique path
in the forward Brownian web starting at 0 at time 0, which is distributed as a standard Brownian motion.

For a fixed point $\omega$ in the probability space of $\bar\W^D$, if
$d(\bar l_{t,\d},\ \B^{0,0}_t)\not\to 0$ as $\d\downarrow 0$, then there exists $\e_\omega>0$
and a sequence $\d_n \downarrow 0$, such that $d(\bar l_{t,\d_n},\ \B^{0,0}_t) >\e_\omega$. In
particular, for $\d_n$ sufficiently small, there exists $(x_n, t_n)\in (\d_n/\s)\Z\times \d_n^2\Z$
with $\sup_n |t_n| <+\infty$,
such that $\bar l_{t_n,\d_n} = l_{t_n, \d_n}=x_n$ and $|x_n-\B^{0,0}_{t_n}| >\e_\omega/2$. Since
$x_n$ is the position of the leftmost zero at time $t_n$ for the rescaled voter model
$\phi^{\Z^-}_{t,\d_n}$, by duality, the backward random walk paths in $\hat\X_{\d_n}$ starting at
$x_n$ and $x_n-\d_n/\s$ at time $t_n$ satisfy $\hat\k^{x_n,t_n}_{0, \d_n}> 0$ and
$\hat\k^{x_n-\d_n/\s,t_n}_{0, \d_n}\leq 0$ at time 0. But by the almost sure convergence of
$\hat\X_{\d_n}$ to
$\hat{\bar\W}$, there exists a subsequence $\d_{n'}$, such that $\hat\k^{x_{n'},t_{n'}}_{t, \d_{n'}}$
and
$\hat\k^{x_{n'}-\d_{n'}/\s,t_{n'}}_{t, \d_{n'}}$ converge respectively to $\hat\k^{x_0,t_0}_t$ and
$\hat\k^{x_0^-,t_0}_t\in\hat{\bar\W}$ for some $(x_0,t_0)$, two paths in $\hat{\bar\W}$ both
starting at $(x_0,t_0)$ with $|t_0|<+\infty$,
$|x_0-\B^{0,0}_{t_0}|>\e_\omega/2$, $\hat\k^{x_0,t_0}_0\geq 0$ and
$\hat\k^{x_0^-,t_0}_0\leq0$. (Note that $x_0\neq \pm \infty$, because for $|t_0|<+\infty$ and
$x_0=\pm\infty$, the only path that can start from $(x_0,t_0)$ is $(f, t_0)$ with $f\equiv x_0$.)
By the non-crossing property of the double Brownian web,
$\hat\k^{x_0,t_0}_t$ and $\hat\k^{x_0^-,t_0}_t$ cannot cross $\B^{0,0}_t$, hence either
$\hat\k^{x_0,t_0}_0=0$ or $\hat\k^{x_0^-,t_0}_0=0$. Therefore
\beqnn
\!\!
\{\bar l_{t,\d}\not\to \B^{0,0}_t \}\subset \{\,\exists\,\hat\k^{x_0,t_0}_t\in\hat{\bar\W}
\hbox{ starting at }(x_0,t_0) \hbox{ with }t_0>0\hbox{ and }\hat\k^{x_0,t_0}_0=0 \}.
\!\!\!\!\!\!\!\!\!\!\!
\eeqnn
The second event in this inclusion has probability zero for the double Brownian web; therefore
almost surely, $\bar l_{t,\d}\to \B^{0,0}_t$ as $\d\to 0$. The same is true for $\bar r_{t,\d}$, and
this proves the theorem for discrete time.

For continuous time, the proof is similar. After applying Skorohod's representation, to prove
$d(\bar l_{t,\d},\ \B^{0,0}_t)\to 0$ almost surely, note that it is sufficient to prove
$d(l_{t,\d},\ \B^{0,0}_t)\to 0$ almost surely instead. The rest of the proof is then the same.

\brm
The conclusions of Theorem \ref{teo:interface} can potentially be used to establish $(T_1)$, the
tightness condition, for $\{\X_\d\}$ and $\{\Xc_\d\}$. In particular, if Cox and Durrett's result
in~\cite{kn:CD} can be extended to establish the conclusions of Theorem \ref{teo:interface} under
the assumption $\E[|Y|^3]<+\infty$, then $\{\X_\d\}$ and $\{\Xc_\d\}$ will be tight, and $\X_\d$
and $\Xc_\d$ will both converge weakly to the Brownian Web under the finite third moment assumption.
\erm

\newpage

\addcontentsline{toc}{chapter}{Bibliography}

%\bibliography{thesis}

\end{document}